%% file: main.tex
\documentclass[fleqn]{article}

\usepackage{arxiv}

\usepackage[utf8]{inputenc} 
\usepackage[T1]{fontenc}    
\usepackage[hidelinks]{hyperref}
\hypersetup{
  linkcolor=[rgb]{0,0,0.5},      
  citecolor=[rgb]{0,0,0.5},      
  urlcolor=[rgb]{0,0,0.5},       
  filecolor=[rgb]{0,0,0.5}
}
\usepackage{url}            
\usepackage{tabularray}
\usepackage{amsfonts}       
\usepackage{nicefrac}       
\usepackage{microtype}      
\usepackage{lipsum}		    
\usepackage{graphicx}
\usepackage[round]{natbib}
\usepackage{subcaption}
\usepackage{enumitem}
\usepackage{amsmath}
\usepackage{bbm}
\usepackage{doi}
\usepackage[ruled,linesnumbered,noend,fillcomment]{algorithm2e} 
\usepackage{setspace}
\usepackage{pdflscape}
\usepackage{multirow}

\title{Drone-Aided Blood Collection Routing Problem: A Column Generation Approach}
\date{}

\author{
Amirhossein Abbaszadeh$^{1,2}$\thanks{Corresponding author: \texttt{amirhossein.abbaszadeh@mail.concordia.ca}}
\And
Hossein Hashemi Doulabi$^{1,2}$
}

\affiliations{
\small
$^{1}$Department of Mechanical, Industrial and Aerospace Engineering, Concordia University, Montréal, Canada\\
$^{2}$Interuniversity Research Center on Enterprise Networks, Logistics and Transportation (CIRRELT), Montréal, Canada
}

\renewcommand{\shorttitle}{Drone-Aided Blood Collection Routing Problem}

\hypersetup{
pdftitle={Drone-Aided Blood Collection Routing Problem},
pdfsubject={Operations Research, Optimization},
pdfauthor={Amirhossein Abbaszadeh, Hossein Hashemi Doulabi},
pdfkeywords={Truck-and-drone routing, Blood collection, Platelet production, Column generation},
}

\begin{document}
\maketitle

\begin{abstract}
	Platelet extraction requires whole blood to be processed within six hours of donation. To meet this deadline, blood collection organizations must optimally route a fleet of vehicles to pick up blood units from donation sites and deliver them to a processing center. This paper introduces a drone-aided blood collection routing problem in which a fleet of trucks, each equipped with a drone, operates in a synchronized manner to collect blood units before their processing time limit expires. Each truck–drone tandem can perform multiple trips throughout the planning horizon, allowing donation sites to be visited repeatedly as new blood units become available over time. We formulate this problem as a mixed-integer linear program that jointly optimizes the routing of trucks and drones, pickup schedules, and timing decisions to maximize the total number of viable blood units collected. We also develop a column generation approach that decomposes the problem into a master problem to select the optimal set of truck–drone tours and a pricing subproblem, which is solved using a tailored memetic algorithm to generate promising new columns. Through a comprehensive computational study, we show the operational benefits of integrating drones into the blood collection system. In addition, we demonstrate the superior performance of the proposed algorithm over Gurobi and two metaheuristics from the literature, namely the hybrid genetic algorithm and the invasive weed optimization, in both the drone-aided and truck-only settings.
\end{abstract}

\keywords{Truck-and-drone routing \and Blood collection \and Platelet production \and Column generation}

\section{Introduction} \label{Sec1}
Once blood is extracted from the human body, it becomes a perishable product that requires careful handling throughout its journey to patients. This journey typically begins at donation sites~(DSs), which are equipped with the facilities required to safely draw blood from donors. After donation, blood is transported to a processing center~(PC) for subsequent testing, processing, and storage. At the PC, blood may either be stored for direct use or separated into its components, such as red blood cells, plasma, and platelets, to address specific clinical needs~\citep{piraban2019survey}.

Among these components, platelets are particularly critical and in high demand, especially for patients undergoing cancer treatments, organ transplants, and major surgeries. Despite their clinical importance, platelets account for less than one percent of total blood volume, and the production of a single therapeutic unit generally requires pooling platelets from four to six blood units~\citep{cotton2015pooled}. In addition, platelet extraction must be completed within six hours of donation to preserve viability. This strict time frame, measured from the moment of donation to the completion of component extraction, is referred to as the processing time limit (PTL)~\citep{ozener2018managing, khameneh2023non}. The time-sensitive nature of this requirement has motivated several studies to propose variants of the vehicle routing problem (VRP) specifically tailored to the blood transportation in support of platelet extraction~\citep[e.g.,][]{yi2003vehicle,doerner2008exact,ozener2018managing}.

Along this line of research, \citet{khameneh2023non} recently studied a multi-trip vehicle routing problem in which a fleet of homogeneous shuttles (trucks) is dispatched to pick up donated blood units from DSs and deliver them to the PC within the PTL, thereby preserving their suitability for platelet extraction. This PTL-constrained pickup-and-delivery problem is hereafter referred to as the \textit{blood collection routing problem} (BCRP). Although the approach adopted by~\citet{khameneh2023non} provides efficient solutions for the BCRP, additional opportunities for improvement exist given the potential of incorporating drones to enhance blood transportation~\citep{callewaert2024limited, yin2024exact, abbaszadeh2025drone}. To make these opportunities explicit, we provide an illustrative example below that motivates our study.

\begin{figure}[h]
    \centering
    \makebox[\textwidth][l]{%
        \hspace{-1.3cm}
        \begin{minipage}{1.12\textwidth}
            \centering

            \begin{subfigure}[b]{0.32\textwidth}
                \centering
                \includegraphics[scale=0.32]{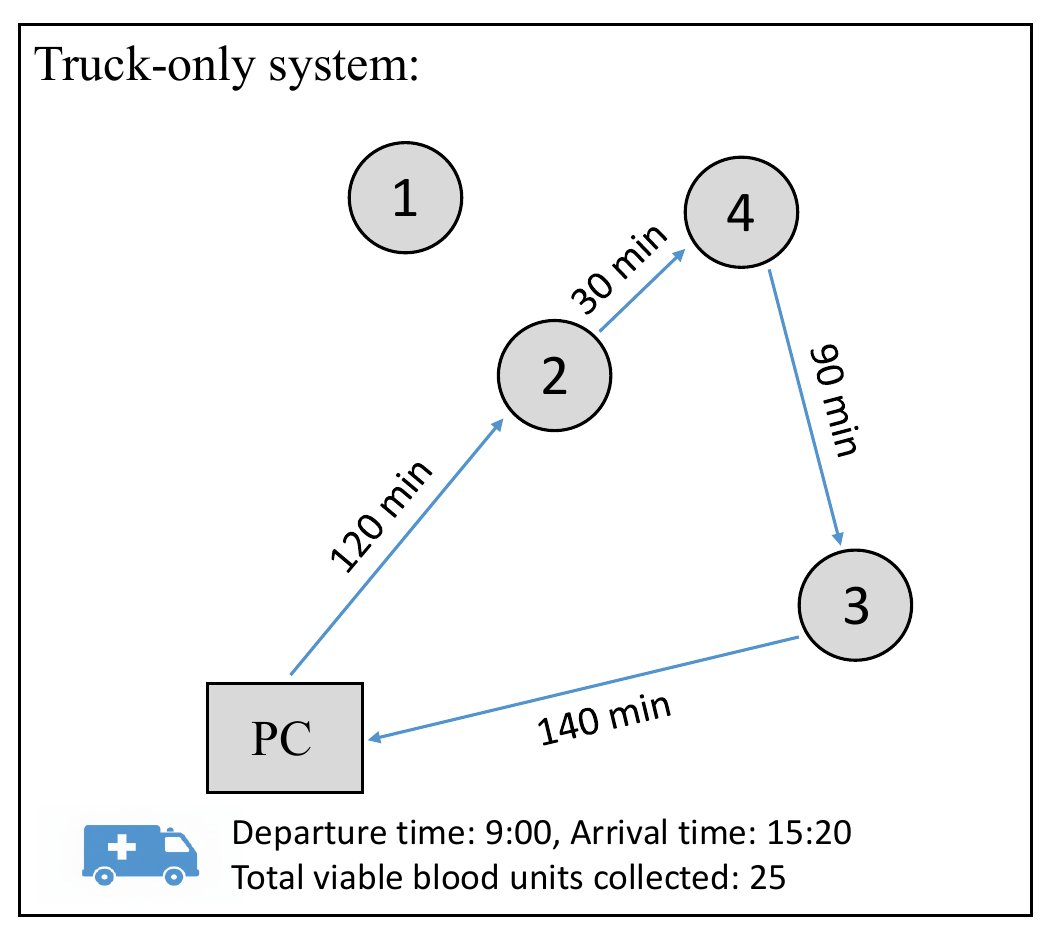}
                \caption{\centering The optimal solution based \\ on \cite{khameneh2023non}}
                \label{fig:Example_PartA}
            \end{subfigure}
            \hfill
            \begin{subfigure}[b]{0.32\textwidth}
                \centering
                \includegraphics[scale=0.32]{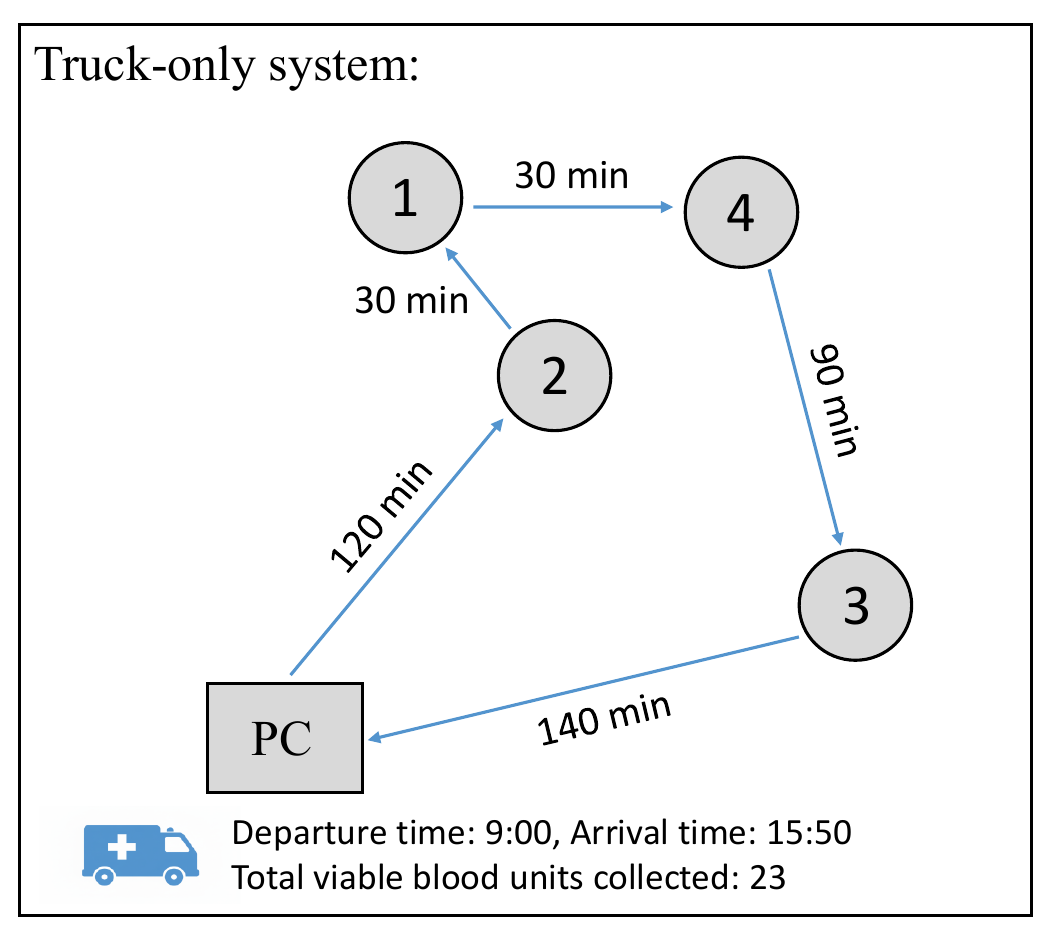}
                \caption{\centering A feasible solution based \\ on \cite{khameneh2023non}}
                \label{fig:Example_PartB}
            \end{subfigure}
            \hfill
            \begin{subfigure}[b]{0.32\textwidth}
                \centering
                \includegraphics[scale=0.32]{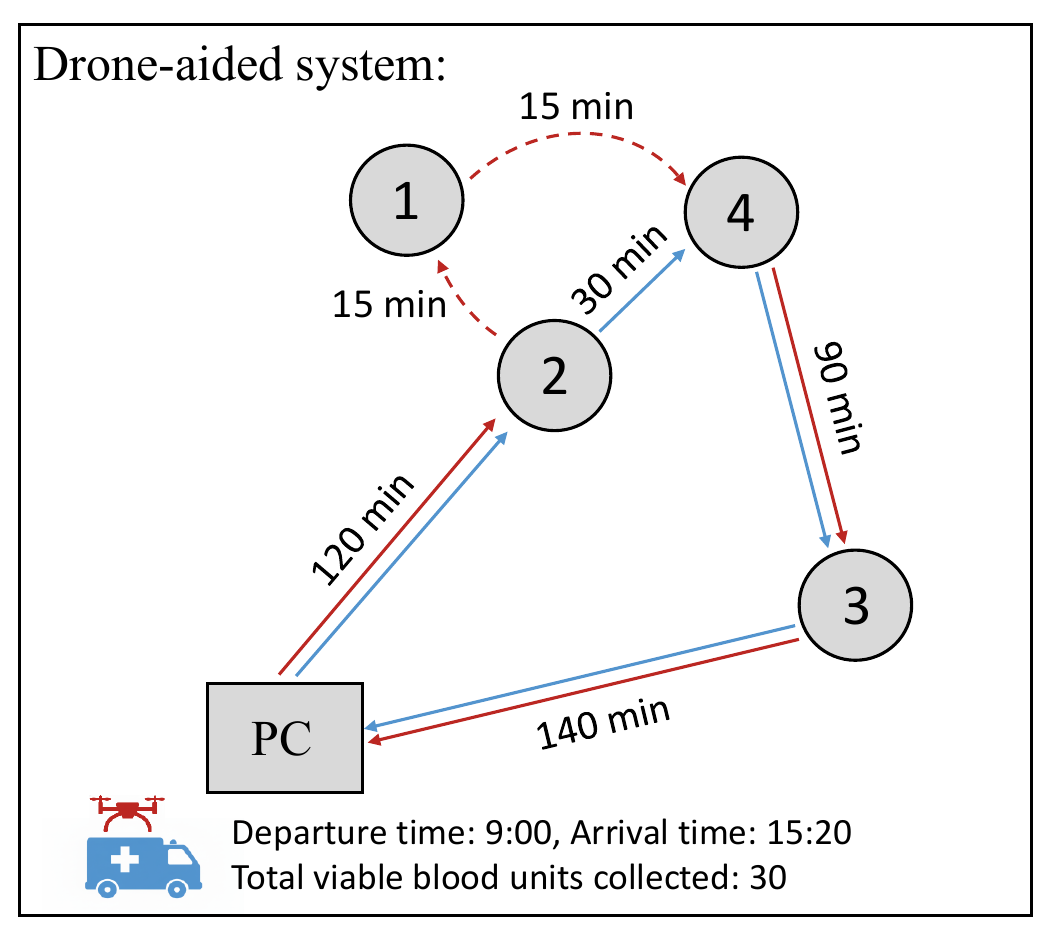}
                \caption{\centering The optimal solution \\ based on our study}
                \label{fig:Example_PartC}
            \end{subfigure}

            \vspace{0.3cm}

            \hspace*{0.5cm}\begin{subfigure}[b]{\textwidth}
                \centering
                \includegraphics[scale=0.45]{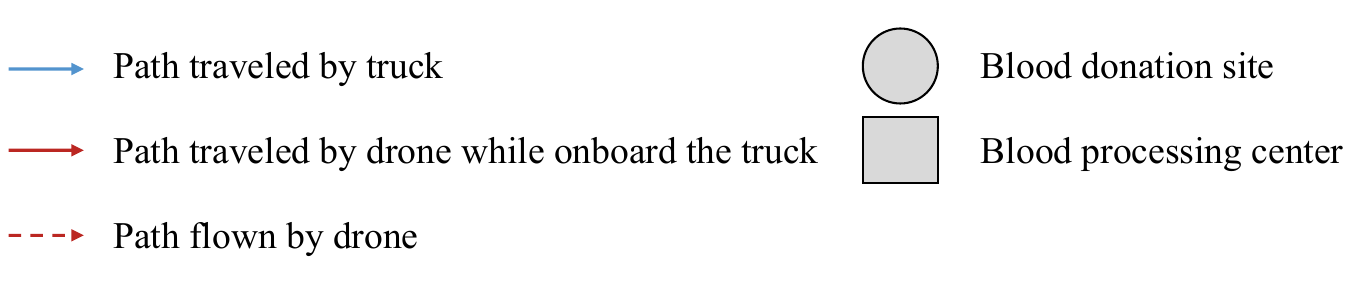}
            \end{subfigure}

            \vspace{0cm}
            \caption{Truck-only and drone-aided solutions for a BCRP instance.}
            \label{fig:Example_ProblemDefinition}
        \end{minipage}
    }
\end{figure}

Consider a simple instance of the BCRP with four DSs and a single PC, with the complete road-travel time matrix provided in~\ref{sec:appTimeMatrix}. Table~\ref{Table:DonationPlan} lists the pre-scheduled blood donation counts at each DS over time, using a 24-hour clock. For example, at 12:30, five blood units are completed and ready for pickup from DS~3. We assume a PTL of five hours, reserving one hour at the PC for platelet extraction, and consider a single truck for blood collection. Fig.~\ref{fig:Example_PartA} illustrates the optimal solution for this instance under the \mbox{truck-only} modeling framework proposed by \citet{khameneh2023non}. In this solution, the truck departs from the PC at~9:00 and visits DSs~2, 4, and~3 at~11:00, 11:30, and 13:00, respectively. During the trip, the truck picks up 6, 6, and 13 blood units from DSs~2, 4, and~3 and returns to the PC at~15:20. All 25 blood units arrive at the PC within the PTL and are therefore suitable for platelet extraction. Hereafter, we refer to such blood units as \textit{viable}.

\input{Tables/Table_Intro}

As observed in the example, and consistent with prior findings in \cite{khameneh2023non}, visiting all donation sites is not necessarily optimal under PTL constraints. In particular, skipping DS~1 is essential to achieving the best solution, as collecting the five blood units completed at 11:00 at this site causes donations at other sites to expire, resulting in an inferior outcome. For example, as illustrated in Fig.~\ref{fig:Example_PartB}, a feasible alternative route that visits DS~2 at 11:00, then DS~1 at 11:30, and subsequently visits DSs~4 and~3 at 12:00 and 13:30 results in only 23 viable blood units upon arrival at the PC at 15:50. Such an inferior outcome occurs because, by the time the truck returns to the PC, the donations completed at 10:30 at DSs~3 and~4 have exceeded the PTL and are no longer suitable for platelet extraction. 

This observation reveals a gap in truck-only operations that can be bridged through the integration of aerial logistics. Specifically, augmenting the collection fleet with a drone allows additional pickups to be performed without delaying the truck’s route. Unlike trucks, drones are not constrained by road congestion or traffic conditions and can travel at higher speeds, enabling parallel and more flexible collection operations. Such a capability is exploited in the drone-aided solution illustrated in Fig.~\ref{fig:Example_PartC}. In the proposed solution, the truck follows the same route and schedule as in the optimal truck-only solution, while a drone is dispatched to visit DS~1 during the truck’s travel between DSs~2 and~4. Assuming that the drone travels at twice the speed of the truck, it reaches DS~1 at 11:15, picks up the five donated blood units, and rejoins the truck at DS~4 at 11:30. The additional pickup at DS~1 does not affect the timing of the truck’s subsequent visits, and all donations picked up during the truck-drone route reach the PC before the PTL expires. Accordingly, the five additional units from DS~1 increase the total number of viable blood units delivered to the PC from~25~to~30, demonstrating the advantage of joint truck–drone operations in the BCRP.

Coordinating truck–drone operations becomes increasingly challenging as the problem size (i.e., the number of DSs and the fleet size) grows. Larger instances require precise synchronization decisions to ensure effective collaboration between trucks and drones under strict PTL constraints. Motivated by these complexities, we address the \textit{drone-aided blood collection routing problem} (DABCRP), and our main contributions are summarized as follows:

\begin{itemize}[label=$\bullet$]
\item To the best of our knowledge, this is the first study to address the synchronization of trucks and drones in the BCRP. The proposed DABCRP extends the truck-only blood collection problem presented in \citet{khameneh2023non} by integrating coordinated drone and truck operations. We show that this integration improves the performance of the collection fleet by increasing the number of blood units delivered to the PC that are suitable for platelet extraction.

\item We present a novel mixed-integer linear programming (MILP) formulation for the DABCRP and develop an efficient column generation approach to solve large-scale instances. The proposed algorithm decomposes the original model into a master problem that selects feasible truck–drone tours and a pricing subproblem that generates promising tours for solution improvement. To overcome the computational difficulty of the pricing subproblem, we design a tailored memetic algorithm that combines evolutionary search with local search operators to efficiently generate high-quality columns.

\item Through a comprehensive computational study, we demonstrate the effectiveness of the proposed solution approach by benchmarking it against the Gurobi solver and two metaheuristic algorithms developed by \citet{khameneh2023non}. The results show that the proposed algorithm not only delivers superior solutions for the drone-aided setting but also consistently improves the best solutions obtained by the competing methods when applied to the truck-only configuration.
\end{itemize}

The remainder of the paper is structured as follows. In Section~\ref{Sec2}, we review the related literature. We formally define the problem in Section~\ref{Sec3} and describe its mathematical formulation in Section~\ref{Sec4}. In Section~\ref{Sec5}, we detail the solution approach. We then evaluate the approach through computational experiments in Section~\ref{Sec6} and conclude the paper in Section~\ref{Sec7} by highlighting the main findings.

\section{Literature review} \label{Sec2}
The current study lies at the intersection of two primary research streams: (i)~routing for blood collection and (ii)~\mbox{drone-aided} routing problems. The following sections provide a detailed review of the relevant literature in each stream.

\subsection{Routing for blood collection}
The blood supply chain consists of multiple stages, including collection, processing, storage, and distribution, each of which requires rigorous planning to ensure the reliable provision of healthcare services. Accordingly, a rich body of literature has addressed strategic and tactical decision-making problems in this context~(see~\citealp{piraban2019survey, agac2024blood, siswanto2025analytical}). Given the focus of our study, in this section we review contributions that include mathematical formulations with routing decisions for various vehicle types to support blood collection operations.

A basic formulation of the BCRP was first introduced by~\cite{prastacos1984blood} as a VRP variant; however, the proposed model fails to account for PTL considerations. \cite{yi2003vehicle} later emphasized the time-sensitive nature of the BCRP through a case study of the American Red Cross, in which all blood deliveries to the PC are subject to a six-hour PTL. In this problem, not all DSs are required to be visited; however, those selected can be served only once within their designated time windows. The author proposed a heuristic algorithm to solve the problem. Motivated by a similar case study conducted with the Austrian Red Cross, \cite{doerner2008exact} advanced this line of research by allowing multiple visits to the same DS. The authors developed several heuristics and a branch-and-bound algorithm to address the problem.  

\cite{mobasher2015coordinating} integrated donor appointment scheduling into the BCRP. In this problem, DSs are grouped into clusters, each served by a single truck, and donated blood must be delivered to the PC within a six-hour PTL to allow platelet extraction. The authors solved the problem using a cluster-first, route-second heuristic. \cite{ozener2018managing} studied a similar problem but considered pre-scheduled donors and addressed clustering and routing decisions in an integrated manner. In their solution approach, the marginal benefit of assigning a DS to a cluster is evaluated based on the resulting increase in platelet production. The authors developed several heuristic algorithms to determine the collection plan after each site–cluster assignment. This research stream was further extended in two subsequent works. First, \cite{khameneh2023non} relaxed the clustering restriction and introduced the non-clustered version of the BCRP studied by \cite{ozener2018managing}, in which trucks are allowed to travel freely among DSs. They solved the problem using two metaheuristic algorithms. Second, \cite{piraban2022multi} extended the framework to account for other types of blood products beyond platelets. To do so, each truck is allowed to perform up to two trips:~one dedicated to collecting units for platelet and cryoprecipitate production, considering an eight-hour PTL, and another to collect any remaining blood that is not subject to PTL constraints. Two metaheuristic algorithms were developed to address the problem. Other studies within this research stream have expanded the problem scope by introducing additional dimensions, although they do not consider the PTL requirement in blood transportation. For example, \cite{mousavi2021designing} developed a sustainable blood supply chain network that incorporates environmental factors and social impacts into the decision-making process. The problem is formulated as a bi-objective stochastic model that minimizes transportation and $\mathrm{CO_2}$ emission costs while maximizing the social benefit obtained from decomposing collected blood into its sub-products. Several multi-objective metaheuristic algorithms were developed to solve the problem. \cite{esmaeili2023exact} integrated facility location, inventory, and vehicle routing decisions into a unified framework. The authors proposed a bi-objective model that simultaneously minimizes total system cost and blood shortages, and is solved using an exact criterion-space search algorithm to generate the complete set of Pareto-optimal solutions.

While all the above-mentioned studies focused on trucks as the primary collection vehicles, other works considered different types of vehicles. In this regard, \citet{csahinyazan2015selective} studied a mobile blood collection problem that employs two vehicle types. The first includes specialized buses, known as \textit{bloodmobiles}, which are equipped with donation beds to facilitate blood donation at public places. The second includes trucks that transport the donated units to the PC. This collection operation allows bloodmobiles to remain at a site for up to three days without risking blood spoilage, as trucks ensure the transfer of donated units at the end of each day. The authors proposed an integer-programming-based heuristic to solve the problem. A similar blood collection system was examined by \cite{rabbani2017solving}, who formulated two mathematical models to support the platelet supply chain. The first model optimizes the locations and schedules of bloodmobiles using a fuzzy multi-objective programming approach. The second is a VRP that optimizes truck tours to deliver the collected blood to the PC before a six-hour PTL expires and is solved using a metaheuristic. \cite{gunpinar2016integer} investigated a bloodmobile routing problem in which the fleet size and the corresponding routes are optimized to satisfy blood demand while minimizing total travel distance. A branch-and-price algorithm was developed to solve the problem. More recently, drones have been incorporated into blood collection operations to improve operational efficiency. Along this line, \cite{kallaj2023integrating} studied a post-disaster blood collection problem in which bloodmobiles visit DSs to collect blood, while drones support the transport of collected units to a crisis-stricken city. The authors formulate the problem as a bi-objective stochastic model that maximizes the total collected blood and minimizes the latest arrival time of vehicles at the crisis region. A multi-objective metaheuristic was used to solve the proposed model. \cite{abbaszadeh2025drone} introduced a synchronized bloodmobile–drone system for mobile blood collection. In this collaborative framework, bloodmobiles conduct on-site collection while drones fly between them to pick up blood bags and deliver them to the PC. This synchronization ensures the timely delivery of blood units for the extraction of time-sensitive components such as platelets. The authors proposed a matheuristic to solve the problem.

Despite extensive research on BCRPs, the synchronized routing of trucks and drones for time-sensitive blood collection operations has not been addressed. This study fills this gap by introducing the first drone-aided BCRP in which trucks and drones are jointly planned to pick up donated blood units from DSs and deliver them to the PC before the PTL expires, thereby ensuring their suitability for platelet extraction. By enabling drones to perform parallel pickups, the proposed collaborative framework expands the operational flexibility of existing routing models and leads to a more effective blood collection system.

\subsection{Drone-aided routing problems}
The rapid evolution of drone technology has opened new opportunities across a wide range of applications. As drones have become more accessible, a significant body of research has emerged focusing on the planning and optimization of drone-aided operations~(see~\citealp{otto2018optimization, telli2023comprehensive}). This interest is especially pronounced in delivery services, where routing decisions play a critical role in achieving operational efficiency~(see~\citealp{shuaibu2025review, zhou2025survey, faramarzzadeh2025literature}).

A seminal contribution in this field was made by \citet{murray2015flying}, who developed two mathematical formulations for drone-aided parcel delivery. The first formulation models the \textit{flying sidekick traveling salesman problem} (FSTSP), in which a truck follows a delivery route while carrying a drone that can be launched to serve individual customers and later rejoin the truck at another location. The second formulation defines the \textit{parallel drone scheduling traveling salesman problem}~(PDSTSP), in which a fleet of drones is dispatched from the depot to serve nearby customers, while a truck independently visits the remaining customers. In a closely related study, \cite{agatz2018optimization} introduced the \textit{traveling salesman problem with drone}~(TSP-D), which differs from the FSTSP by allowing a drone to be launched from and retrieved at the same node. Several subsequent studies advanced the FSTSP, TSP-D, and PDSTSP through strengthened mathematical formulations and improved solution approaches~(e.g.,~\citealp{yurek2018decomposition, dell2020matheuristic, dell2021algorithms, dell2021drone, boccia2021column, roberti2021exact, dell2022exact, blufstein2024decremental, mahmoudinazlou2024hybrid, chen2024learning, zeng2025end, zhang2025benders}). Other studies explored new problem variants under different operational assumptions. For example, \citet{jeong2019truck} extended the FSTSP by incorporating a payload–energy dependency, whereby the weight of a parcel directly affects the drone’s flight duration, while also accounting for time-dependent no-fly zones that restrict flight paths. \citet{nguyen2023parallel} advanced the PDSTSP by introducing collective drones, in which multiple drones can be coupled to jointly serve a single customer, thereby enabling the delivery of heavier or more distant parcels through synchronized drone operations. \citet{hong2025traveling} proposed a TSP-D variant in which the truck can stop at mobile parking locations to enable customer self-pickup, while the drone continues to perform home-delivery services.

Another important research stream emerged with the introduction of the \textit{vehicle routing problem with drones}~(VRP-D) by \citet{wang2017vehicle}. As a fleet-level generalization of the TSP-D, this problem involves a homogeneous fleet of trucks, each carrying a given number of drones, to serve a set of customers. By employing a worst-case analysis, the authors derived theoretical bounds on the potential time savings of \mbox{drone-integrated} fleets compared to traditional truck-only operations. Several follow-up studies refined the VRP-D by developing mathematical formulations and effective solution algorithms~(e.g.,~\citealp{di2017last, schermer2018algorithms, schermer2018variable, schermer2019matheuristic, sacramento2019adaptive, euchi2021hybrid, tamke2021branch, schmidt2025exact}). In parallel, the research stream expanded to include variants reflecting diverse operational settings. For instance, \citet{schermer2019hybrid} introduced a VRP-D with en-route operations, allowing drones to be launched and retrieved at discrete points along arcs. \citet{wang2019vehicle} considered docking hubs as intermediate transfer points that enable drones to be reloaded and reassigned across trucks. \citet{tamke2023vehicle} incorporated speed-dependent energy consumption and treated drone flight speed as a decision variable. \citet{zhou2023exact} formulated a variant in which trucks can stop at selected locations and dispatch their drones to complete several deliveries before continuing onward. \citet{jiang2024multi} proposed a multi-visit flexible-docking VRP-D that allows drones to serve multiple customers per flight, dock with different trucks after each mission, and perform both pickup and delivery operations. \citet{peng2025transportation} studied a time-dependent VRP-D where truck travel speeds fluctuate throughout the day based on traffic conditions, leading to changes in carbon emissions.

Although the drone-aided routing literature is extensive, the problem studied in this paper possesses two main features that existing truck–drone cooperation frameworks are not designed to model. These two features can be characterized as follows: 
\begin{itemize}[label={$\bullet$}]

\item First, pickup operations in the DABCRP are both selective and time-dependent for the following reasons: (i)~ Selectivity arises because, upon arrival at a DS, the vehicle must determine which specific blood donations to pick up among those available at the site. This decision is constrained by the drone’s payload capacity when a drone performs the pickup and, more critically, by the PTL. In particular, since blood units are donated at different times and each has a specific expiration deadline, the vehicle must ensure that all loaded units reach the PC while they are still viable for platelet extraction. (ii)~Time dependency results from the fact that blood donations accumulate over time and are depleted by earlier collections; consequently, the quantity available at any visit depends jointly on the visit time and on all prior pickup decisions at that site.

\item Second, delivery in the DABCRP does not correspond to serving customer locations but to bringing donated blood to the PC before the PTL expires. Accordingly, unlike traditional drone-aided routing models that prioritize travel cost or completion time, the DABCRP prioritizes blood freshness, so that the value of a delivery depends on the time elapsed between donation and arrival at the PC.

\end{itemize}

\section{Problem definition} \label{Sec3}
This section formally defines the DABCRP and is organized into two parts. Section~\ref{subsec:notation} introduces the notation and describes the problem. Section~\ref{subsec:exampleRep} shows how Fig.~\ref{fig:Example_PartC} is expressed mathematically, clarifying the use of the introduced notation.

\subsection{Description of notation} \label{subsec:notation}

Table~\ref{Table_Notation} provides the complete list of sets, parameters, and decision variables. For comparison purposes, the notation is closely aligned with that of \citet{khameneh2023non}, while incorporating the additional elements required to model drone operations. 

\input{Tables/Table_Notation}

The DABCRP is defined on a network consisting of a PC, represented by node~$0$, and a set of DSs denoted by~$\mathcal{V}$. The set of all nodes is given by $\mathcal{V}_0=\mathcal{V}\cup\{0\}$, and the network is connected through the edge set~$E$. Over a single working day defined by the interval $[T_s, T_e]$, blood donations at each DS~$i$ are completed at discrete time slots $t\in T_i=\{a_i,\dots,b_i\}$, where $a_i$ and $b_i$ denote the earliest and latest donation completion times. Blood units remain viable for platelet extraction only if they are delivered to the PC within $K$ time units from their completion time. 

Transport operations are carried out by a fleet of truck–drone tandems denoted by~$F$. Each tandem~$k$ consists of a truck~$k$ and a drone~$k$ that operate in coordination, and their travel times on edge~$(i,j)\in E$ are denoted by $t^v_{ij}$ and $t^d_{ij}$, respectively. A tandem~$k$ can execute a sequence of routes, referred to as a tour, within the planning horizon. Each route~$m \in \mathcal{T}$ represents a complete trip starting and ending at the PC, and the maximum number of routes that a tandem can perform is bounded by $L = \lceil { (T_e - T_s + 2K)}/{\min_{i \in \mathcal{V}} \{ 2t_{0i} \}} \rceil$ \citep{khameneh2023non}.  The binary variable $z_{mk}$ indicates whether tandem $k$ executes route $m$ as part of its tour, while $e_{mk}$ denotes the completion time of that route.

Within route~$m$ of tandem~$k$, each DS~$i$ can be visited at most once, either by the truck~$k$ or by the drone~$k$, and the corresponding departure time from the site is indicated by the variable~$v_{imk}$. Moreover, truck and drone routing decisions for visiting DSs are determined by the $y_{ijmk}$ and $h_{injmk}$ variables, respectively. In the definition of the $h_{injmk}$ variables, a sortie $(i,n,j)$ refers to a triplet of site indices: $i$ is the site from which the drone is launched, $n$ is the site visited for blood pickup, and $j$ is the site where the drone rejoins the truck. The set of feasible sorties is defined as $\mathcal{D} := \big\{(i,n,j) \in \mathcal{V}^3 \mid i, n, j \text{ are distinct}, \, t^d_{in} \leq B, \, t^d_{nj} \leq B \big\}$, which restricts drone operations to those sorties whose individual flight segments can be completed within the battery endurance~$B$. This feasibility definition is supported by the assumption that the drone’s battery is swapped for a fully charged one at each landing in negligible time. 

The objective of the DABCRP is to maximize the total number of viable blood units picked up at DSs and delivered to the PC by the fleet. The parameter $d_{it}$ specifies the number of blood units completed at DS~$i$ at time~$t$. These donations may be picked up either by a truck or a drone, as captured by the decision variables $x_{itmk}$ and $f_{itmk}$, respectively. While truck pickups are not subject to a capacity limitation, the number of units carried by the drone during a sortie is limited by its payload capacity~$Q$. 

\subsection{Notation-based representation of the example} \label{subsec:exampleRep}

In the example presented in Fig.~\ref{fig:Example_PartC}, the working day spans 12 hours, from 8:00 to 20:00, corresponding to a planning horizon of 720 minutes ($T_s = 0$ and $T_e = 720$). The set of DSs is $\mathcal{V}=\{1,2,3,4\}$, which extends to $\mathcal{V}_0 = \mathcal{V} \cup \{0\}$ to account for the PC. There is a path between any pair of nodes in the network; thus, the edge set is given by $E=\{(i,j)\mid i,j\in\mathcal{V}_0,\ i\neq j\}$.

The set of donation times $T_i$ for each DS is determined by converting the clock times reported in Table~\ref{Table:DonationPlan} into elapsed minutes from the start of the working day (i.e., 8:00). Accordingly, a clock time $hh{:}mm$ is mapped to minutes as $t = 60\times(hh - 8) + mm$. For instance, at DS~3, donations are scheduled for 10:30, 11:30, 12:00, and 12:30; hence, $T_3 = \{150, 210, 240, 270\}$, with $a_3 = 150$ and $b_3=270$. The parameter $d_{it}$ is also obtained from Table~\ref{Table:DonationPlan}. For example, DS~3 has $d_{3,150} = 3$, $d_{3,210} = 4$, $d_{3,240} = 1$, and $d_{3,270} = 5$. To ensure viability for platelet extraction, each of these batches must be delivered to the PC within a maximum time limit of $K=300$ minutes from the moment of donation.

A single truck--drone tandem is available for the collection process, i.e., $F=\{1\}$, and its daily tour can include up to $L=6$ routes. Truck travel times $t^v_{ij}$ on each edge $(i,j)\in E$ are given by the road-travel time matrix presented in~\ref{sec:appTimeMatrix}. Drone travel times are computed as $t^d_{ij}=\lceil t^v_{ij}/2 \rceil$, assuming that drones travel twice as fast as trucks on the same path \citep{abbaszadeh2025drone}. This assumption reflects the operational advantage of drones, which are not subject to the road network and traffic constraints faced by ground vehicles. The drone has a payload capacity of $Q=10$ blood units and a maximum battery endurance of $B = 30$ minutes. Under these conditions, the set of feasible drone sorties is $D = \{(1,2,4),(1,4,2),(2,1,4),(2,4,1),(4,1,2),(4,2,1)\}$. 

The optimal solution indicates that the single truck–drone tandem ($k=1$) executes a tour consisting of only one route ($m=1$); thus, $z_{11}=1$. The truck follows the path $0 \rightarrow 2 \rightarrow 4 \rightarrow 3 \rightarrow 0$, while the drone executes sortie $(2,1,4)$ during the truck’s travel from DS~2 to DS~4. This operation is captured by the variables $y_{0211}=y_{2411}=y_{4311}=y_{3011}=1$ and $h_{21411}=1$. The corresponding departure times are $v_{211}=180$, $v_{411}=210$, and $v_{311}=300$ for the truck path, and $v_{111}=195$ for the drone sortie. Finally, the route completion time is $e_{11}=440$. Blood units completed at DSs~2, 4, and~3 are picked up by the truck; therefore, $x_{2,180,1,1} = x_{4,150,1,1} = x_{4,180,1,1} = x_{3,150,1,1} = x_{3,210,1,1} = x_{3,240,1,1} = x_{3,270,1,1} = 1$. Moreover, the blood units completed at DS~1 are picked up by the drone; thus, $f_{1,180,1,1}=1$. This solution results in an objective value of 30, representing the total number of viable blood units delivered to the PC.

\section{Model formulation} \label{Sec4}
We propose the following model for the DABCRP, hereafter referred to as the \textit{drone-aided model} (DAM).
\begingroup
\allowdisplaybreaks
\setlength{\jot}{8pt} 
\begin{align}
&\text{Maximize} \quad \sum_{i \in \mathcal{V}} \sum_{t \in T_{i}} \sum_{m \in \mathcal{T}} \sum_{k \in F} \; \; d_{it} (x_{itmk} +  f_{itmk}) \span \label{eq1} \\
&\sum_{m \in \mathcal{T}} \sum_{k \in F} (x_{itmk} + f_{itmk}) \leq 1, &i \in \mathcal{V}, t \in T_i \label{eq2} \\
&x_{itmk} \leq \sum_{j \in \mathcal{V}_0 \setminus \{i\}} y_{jimk}, & \hspace{-1cm} 
 i\in\mathcal{V}, t\in T_i, m \in \mathcal{T}, k\in F \label{eq3} \\
&f_{itmk} \leq \sum_{j \in \mathcal{V}} \sum_{\substack{n \in \mathcal{V} \\ (j,i,n) \in \mathcal{D}}} h_{jinmk}, &\hspace{0cm} i\in\mathcal{V}, t\in T_i, m \in \mathcal{T}, k\in F  \label{eq4} \\
&\sum_{j \in \mathcal{V}} \sum_{\substack{n \in \mathcal{V} \\ (j,i,n) \in \mathcal{D}}} h_{jinmk} + \sum_{j \in \mathcal{V}_0 \setminus \{i\}} y_{jimk} \leq 1, &i \in \mathcal{V}, m \in \mathcal{T}, k \in F \label{eq5} \\
&\sum_{j \in \mathcal{V}_0 \setminus \{i\}} y_{ijmk} = \sum_{j \in \mathcal{V}_0 \setminus \{i\}} y_{jimk}, &i \in \mathcal{V}, m \in \mathcal{T}, k \in F \label{eq6} \\
& r_{imk} \leq \sum_{j \in \mathcal{V}_0 \setminus \{i\}} y_{jimk}, & i\in \mathcal{V}, m\in \mathcal{T}, k\in F \label{eq7} \\
& \sum_{n \in \mathcal{V}} \sum_{\substack{j \in \mathcal{V}: \\ (i,n,j) \in \mathcal{D}}} h_{injmk}  \leq r_{imk}, & i\in \mathcal{V}, m\in\mathcal{T}, k\in F \label{eq8} \\
& r_{jmk} - r_{imk} + y_{ijmk} \leq  \sum_{n \in \mathcal{V}} \sum_{\substack{n^\prime \in \mathcal{V}: \\ (n,n^\prime,j) \in \mathcal{D}}} h_{n{n^\prime}jmk} - \sum_{n \in \mathcal{V}} \sum_{\substack{n^\prime \in \mathcal{V}: \\ (i,n,n^\prime) \in \mathcal{D}}} h_{inn^\prime mk} + 1, \span \notag \\ 
& \span i,j\in \mathcal{V}: i \neq j, m\in\mathcal{T}, k \in F \label{eq9} \\
& \sum_{j \in \mathcal{V}} \sum_{n \in \mathcal{V}} h_{jnimk} \leq \sum_{j \in \mathcal{V}_0 \setminus\{i\} } y_{jimk}, & i\in \mathcal{V}, m \in \mathcal{T}, k \in F \label{eq10} \\ 
&\sum_{j \in \mathcal{V}} y_{0jmk} = z_{mk}, &m \in \mathcal{T}, k \in F \label{eq11} \\
& \sum_{j \in \mathcal{V}} y_{j0mk} = z_{mk}, &m \in \mathcal{T}, k \in F \label{eq12} \\
&z_{mk} \leq z_{m-1,k}, &m\in\mathcal{T}\setminus\{1\}, k \in F \label{eq13} \\ 
&\sum_{i \in \mathcal{V}_0} \sum_{j \in \mathcal{V}_0 \setminus \{i\}} y_{ijmk} \leq H z_{mk}, &m\in\mathcal{T}, k \in F \label{eq14} \\
& v_{imk} + t^v_{ij} \leq v_{jmk} + H\left(1-y_{ijmk}\right), & \hspace{0cm} i,j\in \mathcal{V}: i \neq j, m\in\mathcal{T}, k \in F \label{eq15} \\
&v_{imk} + t^d_{ij} \leq v_{jmk} + H\Big(1-\sum_{\substack{n \in \mathcal{V}: (i,j,n) \in \mathcal{D}}} h_{ijnmk}\Big),  
& i,j\in \mathcal{V}: i \neq j, m\in\mathcal{T}, k \in F \label{eq16} \\
&v_{imk} + t^d_{ij} \leq v_{jmk} + H\Big(1-\sum_{\substack{n \in \mathcal{V}: (i,j,n) \in \mathcal{D}}} h_{nijmk}\Big),
& i,j\in \mathcal{V}: i \neq j, m\in\mathcal{T}, k \in F \label{eq17} \\
&e_{m-1,k}\mathbbm{1}_{\{m \neq 1\}} + t^{v}_{0i} y_{0imk}  \leq v_{imk}, & i\in\mathcal{V}, m \in \mathcal{T}, k \in F \label{eq18} \\
&v_{imk} + t^{v}_{i0} y_{i0mk} \leq e_{mk}, & i\in\mathcal{V}, m \in \mathcal{T}, k \in F \label{eq19} \\
&e_{mk} \leq t + K + H(1-x_{itmk}-f_{itmk}), &\hspace{-1cm} i\in\mathcal{V}, t\in T_i, m \in \mathcal{T}, k\in F \label{eq20} \\
&v_{imk} \geq t \; (x_{itmk}+f_{itmk}), &\hspace{-1cm} i\in\mathcal{V}, t\in T_i, m \in \mathcal{T}, k\in F \label{eq21} \\
&\sum_{t \in T_i} d_{it} f_{itmk} \leq Q, & i\in\mathcal{V}, m\in\mathcal{T}, k\in F \label{eq22} \\
& x_{itmk}, f_{itmk} \in \{0,1\}, &\hspace{-1cm} i\in \mathcal{V},t\in T_{i}, m\in\mathcal{T}, k\in F \label{eq23} \\
& y_{ijmk} \in \{0,1\}, &\hspace{-2cm} i,j\in \mathcal{V}:i \neq j, m\in\mathcal{T}, k\in F \label{eq24} \\
& h_{injmk} \in \{0,1\}, & \hspace{-3cm} i,j,n \in \mathcal{V}: (i,n,j) \in \mathcal{D}, m\in\mathcal{T}, k\in F \label{eq25} \\
& r_{imk} \in \{0,1\}, & i\in \mathcal{V}, m\in\mathcal{T}, k\in F \label{eq26} \\
& z_{mk} \in \{0,1\}, & m\in\mathcal{T}, k\in F \label{eq27} \\
& v_{imk} \geq 0, & i\in\mathcal{V}, m\in\mathcal{T}, k\in F \label{eq28} \\
& e_{mk} \geq 0, & m\in\mathcal{T}, k\in F \label{eq29}
\end{align}
\endgroup

The objective function~\eqref{eq1} maximizes the total number of donations that \mbox{truck-drone} tandems pick up from the DSs and deliver to the PC within the PTL for platelet extraction. Constraint~\eqref{eq2} requires that the donations completed at a given site and time are picked up at most once. Constraints~\eqref{eq3} and~\eqref{eq4} state that donations can only be picked up from a site if it is visited by a truck or a drone. Constraint~\eqref{eq5} ensures that in each route of a truck-drone tandem, a site can be visited only once. Constraint~\eqref{eq6} enforces flow conservation for each site visited by a truck during a given route. Constraint~\eqref{eq7} states that a site is eligible for drone launch only if the corresponding truck visits that site. Constraint~\eqref{eq8} allows a drone sortie only from a site that is eligible for launch. Constraint~\eqref{eq9} prevents overlapping sorties for the same drone by regulating its onboard status along the associated truck route \citep{dell2022exact}. Specifically, it ensures that once a sortie starts from the launch node, the drone is considered off the truck until it is recovered at the rendezvous node. As a result, no other sortie can be initiated along the truck path between the launch and rendezvous nodes. Constraint~\eqref{eq10} ensures that a drone can be recovered at a site only if the associated truck visits that site. Constraint~\eqref{eq11} and~\eqref{eq12} specify that each active route must form a trip that starts and ends at the PC. Constraint~\eqref{eq13} guarantees that a truck-drone tandem can start a new route only if it has completed the previous one. Constraint~\eqref{eq14} restricts truck movements to active routes by allowing arc traversals only when the corresponding route is performed. Constraints~\eqref{eq15}–\eqref{eq17} enforce time synchronization within each route by linking the truck and drone movement decisions with the corresponding departure times at the donation sites. More precisely, constraint~\eqref{eq15} ensures that if a truck travels directly from one site to another, the departure time from the second site must be greater than or equal to the departure time from the first site plus the truck travel time. Constraint~\eqref{eq16} states that the drone’s flight from the launch node to the service node must precede departure from the service node; hence, the departure time at the service node is no earlier than the departure time at the launch node plus the drone travel time. Constraint~\eqref{eq17} ensures that the truck can only depart from the rendezvous node after the drone has returned from the service node. Accordingly, the departure time at the rendezvous node must be greater than or equal to the departure time at the service node plus the drone travel time. Constraint~\eqref{eq18} imposes that the departure time at the first visited site in each route accounts for the travel from the PC and, when applicable, the completion of the previous route. In this constraint, $\mathbbm{1}_{\{m \neq 1\}}$ is an indicator that equals~1 when $m \neq 1$ and~0 otherwise. Constraint~\eqref{eq19} ensures that the completion time of a route reflects the departure time from the last visited site plus the travel back to the PC. Constraint~\eqref{eq20} enforces the platelet viability limit by requiring that if donations are picked up at a site, the corresponding truck-drone tandem must return to the PC within $K$ time units of the donation time. Constraint~\eqref{eq21} specifies that if donations are picked up at a site, the departure from that site must occur no earlier than the corresponding donation time. Constraint~\eqref{eq22} ensures that the number of donations picked up by a drone during a sortie does not exceed its payload capacity. Constraints~\eqref{eq23}-\eqref{eq29} indicate the domains of the decision variables.

\section{Column generation algorithm} \label{Sec5}
The DAM is computationally intractable for generic solvers due to the extremely large number of decision variables and the complex synchronization constraints. To address this challenge, we develop a column generation algorithm~(CGA) that decomposes the DAM into a master problem~(MP) and a pricing subproblem~(PSP). The former is a binary integer program that selects tours, whereas the latter is a mixed-integer program that identifies new tours promising for solution improvement. 

In the following, we first present the formulations of the MP and the PSP. We then introduce a tailored memetic algorithm designed to solve the PSP efficiently. Finally, we describe the procedure for obtaining the final integer solution. To ensure clarity, we reiterate that a route refers to a single trip starting and ending at the PC, whereas a tour represents the complete sequence of routes performed by a tandem over the planning horizon.

\subsection{Master problem formulation}
Let $\mathcal{P}$ denote the set of all feasible tours (columns) that can be operated by 
a \mbox{truck-drone} tandem. The column vector $\mathbf{C}_p$, representing the $p$th tour, is defined as follows:

\begin{align}
\mathbf{C}_p := \left( u^{it}_p \right)_{i \in \mathcal{V}, t \in T_i}^\top,
\end{align}
where $u^{it}_p$ is a binary parameter equal to 1 if the donations completed at DS~$i$ at time~$t$  are picked up in tour~$p$, and 0 otherwise. We formulate the MP as follows:

\begingroup
\allowdisplaybreaks
\setlength{\jot}{8pt} 
\begin{align}
\text{MP:} \quad &\text{Maximize} \quad \sum_{i \in \mathcal{V}} \sum_{t \in T_i} d_{it} \delta_{it} \span \label{eq31} \\
&\text{subject to:} \ \; \delta_{it} \leq \sum_{p \in \mathcal{P}} u^{it}_{p} \lambda_{p}, & i \in \mathcal{V},\ t \in T_i \label{eq32} \\
&\sum_{p \in \mathcal{P}} \lambda_{p} \leq |F|, \span \label{eq33}  \\
&\delta_{it} \in \{0,1\}, & i \in \mathcal{V},\ t \in T_i \label{eq34} \\
&\lambda_{p} \in \{0,1\}, & p \in \mathcal{P} \label{eq35}
\end{align}
\endgroup

Here, $\delta_{it}$ is a binary variable indicating whether the donations completed at DS~$i$ at time~$t$  are collected~(i.e., picked up by a vehicle and delivered to the PC), and $\lambda_{p}$ is a binary variable representing whether tour~$p$ is selected. The objective function~\eqref{eq31} maximizes the total number of blood units collected. Constraint~\eqref{eq32} ensures that donations can only be marked as collected if they are covered by at least one of the selected tours. Constraint~\eqref{eq33} restricts the number of selected tours to the fleet size, ensuring feasibility with respect to the number of available truck–drone tandems. Finally, constraints~\eqref{eq34} and~\eqref{eq35} enforce the binary nature of the decision variables.

Despite its compact form, the MP involves a large number of potential tours even for small instances, making complete enumeration computationally intractable. However, by relaxing the integrality constraints on the variables, the resulting linear master problem~(LMP) can be efficiently solved using column generation. The column generation starts from a restricted version of the problem, called the restricted linear master problem~(RLMP), which contains only a limited subset of tours $\mathcal{P^\prime} \subseteq \mathcal{P}$. The RLMP is solved to optimality, the corresponding dual values are obtained, and these values are then used in the PSP to search for new promising tours that have strictly positive reduced cost. If such tours are identified, they are added to $\mathcal{P^\prime}$ and the RLMP is re-optimized. This process is repeated until no additional promising tours can be generated, at which point the solution of the RLMP is optimal for the LMP~\citep{desaulniers2006column}.

\subsection{Pricing subproblem formulation}

Let \( \pi_{it} \in \mathbb{R}^{+}_{0} \) and \( \mu \in \mathbb{R}^{+}_{0} \) denote the dual values associated with constraints \eqref{eq32} and \eqref{eq33}, respectively. We formulate the PSP as follows:
\begingroup
\allowdisplaybreaks
\setlength{\jot}{8pt} 
\begin{align}
\text{PSP $(\pi, \mu, k)$:} \; &\text{Maximize} \; \left\{ \sum_{i \in \mathcal{V}} \sum_{t \in T_i} \sum_{m \in \mathcal{T}} \pi_{it} (x^{it}_{mk} + f^{it}_{mk}) - \mu \right\} \span \label{eq36} \\ 
&\text{subject to:} \ \; \text{Constraints } \eqref{eq2} - \eqref{eq28}, \text{ where } F \text{ is restricted to } \{k\}
\end{align}
\endgroup

The PSP seeks to identify a tour that can potentially improve the current solution of the RLMP. The objective function maximizes the reduced cost of a new tour, which is driven by the dual values $\pi_{it}$ and $\mu$. The term $\pi_{it}$ reflects the marginal benefit of covering donation $d_{it}$ under the current LP solution, while $\mu$ penalizes the use of an additional truck-drone tandem. A tour with strictly positive reduced cost corresponds to a feasible column that is worth adding to the RLMP. If no such column exists, the current LP solution is optimal for the LMP.


\subsection{Memetic algorithm for solving the pricing subproblem}
Although the PSP optimizes the tour of a single truck–drone tandem, it remains computationally demanding for generic solvers and constitutes the main bottleneck of the CGA. To overcome this challenge, we develop a memetic algorithm (MA) to efficiently generate high-quality solutions for the PSP within a reasonable computational time. The proposed MA builds upon the genetic algorithm (GA) developed by~\cite{khameneh2023non} through the integration of local search (LS) operators introduced by~\cite{ha2020hybrid}.

The MA evolves a population of feasible tours operated by a single truck–drone tandem over a finite number of generations, defined by the parameter \texttt{genLimit}. Each generation adopts one of two search strategies: (i)~solution improvement through crossover and education, or (ii)~population diversification through mutation. In the first strategy, route segments from two parent tours are recombined to generate new offspring, and a set of LS operators is subsequently applied to refine the internal structure of each route. In the second strategy, a mutation operator is triggered to replace a fraction of the least-fit tours with newly generated random ones. This balance between exploitation and exploration enables the MA to efficiently search for high-quality tours within the solution space. The overall framework of the proposed MA is summarized in Algorithm~\ref{alg:Memetic}, with detailed descriptions of its main components presented in Sections~\ref{sec: MA_SlnEncoding} to \ref{sec: MA_Mutation}.

\begin{center}
\begin{minipage}{0.97\textwidth}
\begin{algorithm}[H]
    \small
    \caption{Memetic Algorithm for Solving the PSP}
    \label{alg:Memetic}
    \setstretch{1.30}
    \DontPrintSemicolon

    \KwIn{\textit{genLimit}, \textit{popSize}, \textit{muRate}, \textit{muVolume}}
    
    \KwOut{Best tour found}

    Generate the route pool \;
    $P \gets \emptyset$\;
    \For{$\textit{cnt} = 1$ \KwTo \textit{popSize}}{
        $C \gets$ Construct a random tour using the route pool\;
        $C \gets$ Apply TRP on $C$  \hspace{3.9cm}\tcp{Timing Refinement Procedure (TRP)}
        $P \gets P \cup \{C\}$\;
    }
  
    \For{$\textit{gen} = 1$ \KwTo \textit{genLimit}}{
        $\textit{rnd} \gets \text{random} \; (0, 1)$ \;
        \If{$\textit{rnd} > \textit{muRate}$}{
            $P_1, P_2 \gets$ \textit{RouletteWheelSelection} ($P$)\;
            $O_1, O_2 \gets$ \textit{OnePointCrossover} ($P_1, P_2$)\;
            \ForEach{offspring $O \in \{O_1, O_2\}$}{
                $O \gets$ Apply RARP on $O$ \hspace{2.5cm}\tcp{Route-Alignment Repair Procedure (RARP)}
                \ForEach{route $R \in O$}{
                    $R \gets$ Educate $R$ \hspace{3cm}\tcp{Local Search Using $M_1$ to $M_{16}$}
                }
                $O \gets$ Apply TRP on $O$ \;
                Add $O$ to $P$ and remove the least-fit tour\;
            }
        }
        \Else{
            $P_{weak} \gets$ Select \textit{muVolume} least-fit tours from $P$\;
            \ForEach{$C \in P_{weak}$}{
                $C_{new} \gets$ Construct a random tour using the route pool\;
                $C_{new} \gets$ Apply TRP on $C_{new}$\;
                Replace $C$ with $C_{new}$ in $P$\;
            }
        }
        $C_{best} \gets \text{argmax}_{C \in P} \{ \text{fitness} \; (C) \}$\;
    }
    \Return $C_{best}$\;
\end{algorithm}
\end{minipage}
\end{center}

\subsubsection{Solution encoding} \label{sec: MA_SlnEncoding}
\begin{figure}[h]
    \centering
    \includegraphics[width=\linewidth]{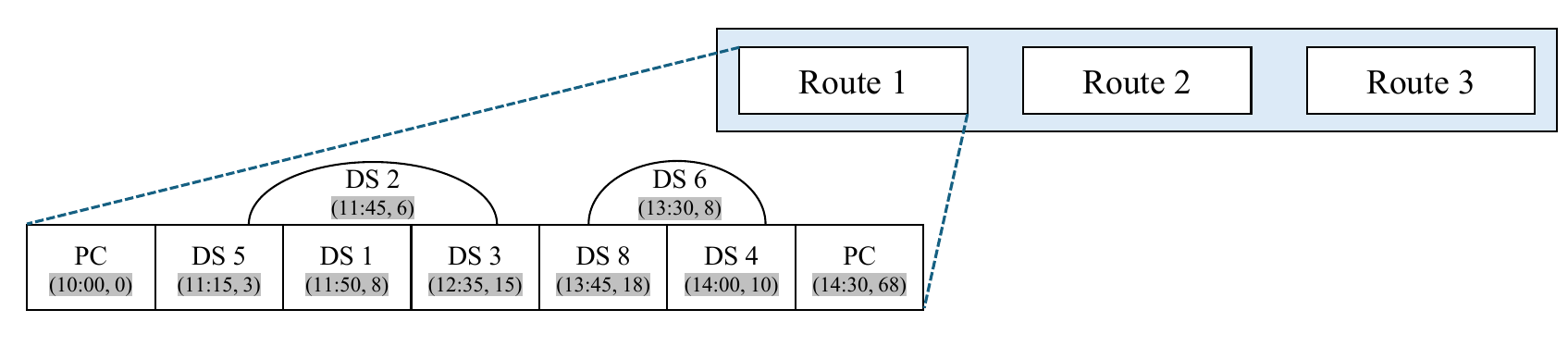}
    \captionsetup{justification=centering} 
    \caption{A sample chromosome encoding.}
    \label{fig: slnEncoding}
\end{figure}

Similar to the GA, the MA encodes each solution as a chromosome. Each chromosome corresponds to a tour operated by a single truck–drone tandem, while each gene represents an individual route performed within that tour. Fig.~\ref{fig: slnEncoding} illustrates an example of a chromosome adopted in the proposed MA. The chromosome consists of three genes (i.e., Route~1, Route~2, and Route~3), each defining a feasible route. Each route represents a coordinated truck–drone operation, including truck stops and drone sorties. In the figure, truck-served DSs are represented by rectangular nodes, whereas drone-served DSs are depicted with curved shapes. Moreover, each DS name is followed by a gray-highlighted tuple that contains two elements: the departure time of the tandem from that site and the amount of blood picked up at that location. For the PC, the tuple indicates the completion time of the route and the total amount of blood delivered to the processing center. In the remainder of this paper, we refer to this pair of values as the \textit{status tuple}. Note that, in the status tuple of a DS, the picked-up quantity is determined based on the departure time from the site. When a capacity limitation applies, such as for drone-served sites, fresher blood units are prioritized so that the most recent donations within the payload limit are picked up.

As depicted in Fig.~\ref{fig: slnEncoding}, in Route~1, the truck visits DSs~5,~1, and~3 at 11:15, 11:50, and 12:35, picking up 3, 8, and 15 blood units, respectively. Meanwhile, the drone performs a sortie: it is launched from DS~5 at 11:15, visits DS~2 at 11:45 to pick up 6 blood units, and rejoins the truck at DS~3 at 12:35. The trip then continues as the truck visits DS~8 and~4, while the drone performs another sortie to pick up blood from DS~6. Route~1 is completed at 14:30, with a total of 68 viable blood units delivered to the PC. The truck–drone tandem then proceeds to perform Route~2 and Route~3. 

\subsubsection{Initial population}
The MA begins by constructing an initial population of feasible tours, whose cardinality is defined by the parameter \texttt{popSize}. The initial population is generated by sequentially combining routes that are randomly selected from a pool of feasible routes. This pool, referred to as the \textit{route pool}, is constructed through a progressive expansion procedure. This procedure first creates all single-visit, truck-only routes, each representing a direct road trip from the PC to a single DS and back to the PC. These elementary routes serve as the initial seeds of the pool. 

The pool is then iteratively expanded through two types of extensions. The first is the \textit{road-travel extension}, in which each existing route in the pool is lengthened by inserting an unvisited DS before the return to the PC. The second is the \textit{aerial-travel extension}, in which a drone sortie is incorporated into each existing route by selecting feasible launch and rendezvous nodes along the path and assigning the drone to visit an unserved DS between them. In this phase, the number of intermediate stops that the truck can make while the drone is performing its sortie is restricted by the parameter \texttt{divLimit}. This diversion limit promotes efficient drone utilization by allowing the drone to return earlier and perform additional sorties within the same tour. The expansion process continues until the number of generated routes reaches the predefined parameter \texttt{poolSize}.

Once the route pool is constructed, a tour is generated as follows. First, a route is randomly selected from the pool, and the corresponding status tuples are computed under the assumption that this is the only route in the tour and that the tour starts at time zero. Next, additional routes are randomly drawn from the pool and appended sequentially under the assumption that no idle time exists between consecutive routes. After each addition, the status tuples of the DSs in the newly added route are computed, and the completion time of the tour (i.e., the completion time of its last route) is updated accordingly. The addition of routes continues until the completion time of the tour reaches the upper limit $FT$. This limit is defined as $FT = LDCT + K$, where $LDCT$ denotes the latest donation completion time within the planning horizon. Any route whose completion time exceeds $FT$ would deliver blood units that are no longer suitable for platelet extraction, as donations completed at $LDCT$ would already exceed the processing time limit $K$. Finally, the fitness of the constructed tour is calculated according to equation~\eqref{eq36}.

Before being added to the initial population, each constructed tour undergoes a \textit{timing refinement procedure} (TRP) aimed at improving its fitness. In this procedure, the idle time of the truck–drone tandem performing the tour is calculated as $IT=FT-CTT$,  where $CTT$ denotes the completion time of the tour. When $IT > 0$, the available time is distributed among the routes using the Dirichlet distribution \citep{lin2016dirichlet}. This operation is repeated \texttt{refIter} times to generate alternative timing patterns and retain the best configuration. The refined tour with the highest fitness value is then added to the initial population.

\subsubsection{Crossover operator}
The MA employs a one-point crossover operator that generates new candidate solutions by exchanging routes between two existing tours selected from the current population. In each generation, the crossover operator is invoked with a probability of $1 - \texttt{muRate}$. When invoked, two tours are selected as parents using the roulette-wheel selection method, where the probability of selection for each tour is proportional to its fitness. Once the parents are selected, the one-point crossover is performed as follows. For each parent, a random cut point is determined along its sequence of routes, and two new tours are produced by exchanging the corresponding route segments. The first offspring inherits the initial set of routes from the first parent and the remaining routes from the second, while the second offspring receives the complementary combination.

After recombination, the resulting offspring may become infeasible in terms of temporal continuity between consecutive routes or in the cumulative amount of blood collected. To restore feasibility, a \textit{route-alignment repair procedure} (RARP) is applied. This procedure first recalculates the status tuples of nodes under the assumption that no idle time exists between consecutive routes and that the tour begins at time zero. It then removes any route whose completion time exceeds the upper limit $FT$ and finally applies the TRP to refine the time-related variables of the offspring. 

\subsubsection{Education phase}
Following the crossover, each offspring undergoes an education phase aimed at improving its fitness value. This phase is designed to refine the internal structure of routes (i.e., the sequencing of DSs within routes), which remain unchanged after executing the crossover operation. During the education phase, all routes of the offspring are sequentially processed. For each route, sixteen LS operators (moves) adapted from \citet{ha2020hybrid} are applied to refine its structure and enhance the overall quality of the offspring. These operators, denoted M$_1$ to M$_{16}$, are described in detail in~\ref{sec:appLSOperators}.

The education phase starts by iteratively applying the LS operators to the first route of the offspring. At each iteration, an active operator M$_{\ell}$ is selected with probability $w_{\ell}/W$, where $w_{\ell}$ denotes its weight and $W = \sum_{\ell=1}^{16} w_{\ell}$. If the selected operator fails to yield an improvement, the solution is deemed locally optimal within the corresponding neighborhood, and the operator is deactivated and excluded from further selection. Conversely, if the selected operator results in an improvement, the previously deactivated operators are again made available for selection. The process continues until all operators become inactive, after which the current route is considered locally optimized. The same procedure is then repeated sequentially for the remaining routes of the offspring.

The operator weights are dynamically adjusted using statistics collected by a short-term memory \citep[see][]{cattaruzza2014memetic}. For each operator M$_{\ell}$, the memory maintains two counters: (i)~$a_{\ell}$, representing the number of times the operator has been applied, and (ii)~$s_{\ell}$, denoting the number of successful applications that improved the solution. After every $\omega$ executed moves ($\omega=100$ in this study), the operator weights are updated according to the rule $w_{\ell} \leftarrow w_{\ell} + s_{\ell}/a_{\ell}$. Once the update is performed, the attempt and success counters are reset to zero for the next evaluation period. This adaptive mechanism gradually reinforces operators that frequently produce improving moves, while maintaining diversity by keeping all operators eligible for selection.


After all routes of the offspring have been educated, the TRP is applied once more to refine the time-related variables. The resulting educated offspring is added to the population, and the least-fit tour is eliminated.

\subsubsection{Mutation operator} \label{sec: MA_Mutation}
To maintain population diversity and prevent premature convergence, the MA incorporates a mutation operator that occasionally replaces part of the population with newly generated solutions. This operator promotes diversification at the tour level and is activated in each generation with a probability of \texttt{muRate}. Once activated, a portion of the least-fit tours, determined by the parameter \texttt{muVolume}, is removed from the population and replaced with newly generated tours. Each new tour is generated in the same manner as during population initialization, where a random sequence of routes from the pre-generated route pool is concatenated, and the TRP is subsequently executed to refine the time-related variables.

\subsection{Finding the final integer solution}
Upon convergence of the column generation procedure, the solution of the RLMP is examined for integrality. If it satisfies the integrality constraints~\eqref{eq34} and~\eqref{eq35}, it constitutes the final solution to the DAM. Otherwise, when the RLMP yields a fractional solution, a post-optimization phase is performed to obtain a feasible integer solution. In this phase, the integrality constraints are reinstated in the RLMP, and the resulting problem is solved to optimality using a commercial solver. This problem determines the combination of generated tours that maximizes the total number of viable blood units collected while satisfying constraints~\eqref{eq32} to~\eqref{eq35}. The resulting solution serves as the final output of the CGA algorithm.

\section{Computational study} \label{Sec6}

This section presents a comprehensive computational study designed to evaluate both the performance of the CGA and the operational impact of integrating drones into the blood collection routing problem. To this end, we first describe the generation of a diverse set of benchmark instances and the calibration of algorithmic parameters. We then compare the proposed CGA against two metaheuristics from the literature and a commercial solver under both truck-only and drone-aided systems. Finally, a sensitivity analysis is conducted to examine the effect of fleet size and drone capabilities on overall blood collection performance. All algorithms and models were implemented in C++, and each experiment was performed on a machine equipped with two Intel Xeon Gold 6148 processors at 2.40~GHz and 50~GB of RAM, using a single core. Moreover, Gurobi~12.0 was used as the solver.

\subsection{Generated instances}
To evaluate the performance of the proposed solution approach, a comprehensive set of test instances was generated following a procedure similar to the one proposed by \citet{ozener2018managing} and later adopted by \citet{khameneh2023non}. The test set comprises three instance classes (C1 to C3), each corresponding to a different problem size. Each class contains ten instances, and each instance is characterized by two key features: the number of DSs and the number of available truck–drone tandems. In class C1, the smallest problems are considered, ranging from 3 to 21 DSs and 1 to 4 tandems. Class C2 introduces medium-sized configurations with 23 to 41 DSs and 5 to 8 tandems. Finally, class C3 represents the largest problems, containing between 43 and 61 DSs and 9 to 12 tandems.

Each instance is defined on a~$1000\times1000$~grid representing the service area, which contains a single PC and the corresponding set of DSs. The DS coordinates are generated in a manner that reflects realistic geographic patterns with both concentrated and scattered locations. To achieve this, two coordinate-generation methods are employed: the first produces clustered sites and the second dispersed ones, with half of the DSs generated by each method. In the first method, each DS is randomly assigned to one of four predefined zones whose centers are located at the midpoints of the grid’s quadrants, representing densely populated regions distributed across the service area. Then, the position of the DS within the selected zone is determined in polar form, where the angular component $\theta$ is drawn from a uniform distribution over~$[0,2\pi]$ and the radial distance $\rho$ is sampled from a normal distribution with mean~0 and standard deviation~100. Finally, the Cartesian coordinates of the site are computed as $(x, y) = (x_c + \rho\cos(\theta), y_c + \rho\sin(\theta))$, where $(x_c, y_c)$ denotes the center of the chosen zone. In the second method, both the $x$- and $y$-coordinates are drawn independently from uniform distributions on $[0,1000]$, resulting in a dispersed set of sites that emulate remote regions within the service area.

The travel distances between nodes are computed using the Euclidean metric. The road travel times are then derived by scaling these distances so that a truck can complete two round-trips between the PC and the farthest DS within the planning horizon. Aerial travel times are obtained by applying a speed-scaling factor~$\alpha=2$, implying that drones are twice as fast as trucks over the same distance \citep{roberti2021exact, abbaszadeh2025drone}. This reflects their advantage of operating independently from road networks and traffic conditions. Each drone has a maximum flight time of 30 minutes and a payload capacity of 5 kg, consistent with the specifications of the commercially available T-DRONES MX860. Considering that a standard blood bag contains approximately 450~mL and weighs about 477~g (based on a density of 1.06~g/mL), each drone can carry up to 10 blood bags per trip.

The planning horizon covers an 800-minute workday (approximately 13~hours), with the start and end times set to $T_s=0$ and $T_e=800$, respectively. The processing time limit for platelet extraction is set to $K=300$~minutes. For each DS~$i$, the earliest and latest donation completion times, $a_i$ and $b_i$, are generated from uniform distributions on $[T_s,T_e/3]$ and $[2T_e/3,T_e]$, respectively, ensuring temporal diversity among sites. To represent heterogeneity in daily donation activity, each DS is randomly classified as having either a low or high donation rate, with equal probability. Low-activity sites have the number of daily donations drawn uniformly from $[1,25]$, whereas high-activity sites have this value drawn from $[25,50]$. The donation completion times within each site are distributed according to two temporal patterns: (i) for half of the DSs, donation times are uniformly generated between $a_i$ and $b_i$; and (ii) for the remaining half, donation times follow a normal distribution centered at $(a_i+b_i)/2$. This design introduces both uniform and peak-oriented temporal profiles, resulting in a diverse set of temporal conditions for algorithmic evaluation.

\subsection{Parameter tuning}
The performance of the CGA depends on several algorithmic parameters. To identify an effective configuration, a Taguchi experimental design was employed using an L27 orthogonal array, which enables efficient exploration of multi-factor interactions with a limited number of experiments. The parameters were tuned through a three-level experimental design, with factor levels summarized in Table~\ref{tab:Taguchi}.

\input{Tables/Table_Taguchi}

\begin{figure}[h]
    \centering
    \includegraphics[scale=0.90]{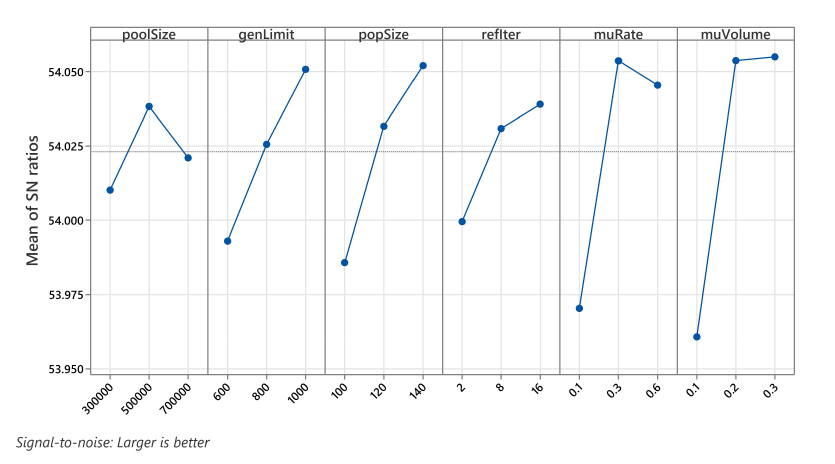}
    \captionsetup{justification=centering} 
    \caption{S/N ratio chart.}
    \label{fig:SN_Ratio}
\end{figure}

The experiments were conducted on the sixth instance of class C2, which represents a medium-sized problem. For each configuration, the CGA was executed ten times, and the objective value was used to compute the signal-to-noise (S/N) ratio in Minitab 22.4.0. As shown in Fig.~\ref{fig:SN_Ratio}, higher values of \texttt{genLimit}, \texttt{popSize}, \texttt{refIter}, and \texttt{muVolume} improve performance, whereas moderate levels of \texttt{poolSize} and \texttt{muRate}  maintain an effective search balance.

\subsection{Benchmark algorithms}
To evaluate the effectiveness of the CGA, we compare its performance with two metaheuristics introduced by \citet{khameneh2023non}, namely the hybrid genetic algorithm (HGA) and the invasive weed optimization (IWO). In what follows, we first outline their main mechanisms and then explain how they are tailored to the DABCRP.

\begin{itemize}
   \item \textbf{HGA} constructs a multi-truck solution through a sequential process, with each iteration focusing on a single truck. For each truck, the algorithm forms the initial population of chromosomes (tours) by randomly combining routes selected from a pre-generated route pool, using a procedure similar to that applied in our study. The fitness of each chromosome is evaluated according to the total number of viable blood donations picked up and delivered within the PTL. The algorithm randomly applies either crossover or mutation for a predefined number of generations. When crossover is performed, route segments between two parent tours are exchanged, whereas the mutation operator replaces a subset of tours with newly generated ones formed by combining routes drawn from the route pool. Once the evolutionary process is completed, a simulated annealing algorithm is applied to refine the best tour by optimizing the departure times of routes and further improving solution quality. After determining the tour for one truck, the set of available blood donations at each site is updated to account for the collected units. The procedure continues with the next truck and is repeated iteratively until all vehicles have been assigned tours.

   \item \textbf{IWO} employs the same sequential procedure as the HGA to construct tours for the fleet of trucks. The algorithm is inspired by the colonization and reproduction mechanisms of invasive weeds, where each solution for a single truck, referred to as a \textit{weed}, represents a one-day tour. The algorithm begins with an initial population of weeds generated from the same pre-defined route pool and evolves iteratively through three key operators: reproduction, spatial dispersal, and competitive exclusion. During reproduction, each weed produces a number of offspring, referred to as \textit{seeds}, in proportion to its fitness. 
   Each seed is generated by either (i) swapping the positions of two routes within the tour or (ii) occasionally replacing an entire route with another drawn from the route pool. The spatial dispersal operator introduces variation by perturbing departure times, with the magnitude of these perturbations gradually decreasing over time to promote convergence. Finally, the competitive exclusion operator ensures that only the fittest individuals survive to the next generation, thereby maintaining a constant population size. Once the best tour for a truck is identified, the set of available blood donations is updated to remove the collected units, and the algorithm proceeds to the next truck. This process is repeated iteratively until all vehicles in the fleet have been assigned tours.
   
\end{itemize}

Both benchmark algorithms are originally designed for truck-only blood collection. To make them capable of solving the proposed drone-aided problem and enable a meaningful comparison with the CGA, their solution representations are adapted while preserving their core search logic. Each route in a chromosome (for HGA) or weed (for IWO) is redefined as a coordinated truck–drone route consisting of truck stops and feasible drone sorties. Accordingly, each tour represents a one-day operation of a single truck–drone tandem, encoded as described in Section~\ref{sec: MA_SlnEncoding}. All other algorithmic operators remain unchanged from \citet{khameneh2023non} to ensure consistency in performance evaluation, and the required parameters are tuned using the Taguchi method to identify an effective configuration.

It is worth noting that our proposed CGA can also be adapted to the truck-only blood collection setting. To achieve this, two minor adjustments are required. First, during the initialization phase, the aerial-travel extension is disabled when constructing the route pool, ensuring that no drone sorties are included in the generated routes. Second, in the education phase, the LS operators dedicated to drone-related adjustments (i.e., operators M$_9$-M$_{16}$) are deactivated, preventing the introduction of drone operations during local search. With these modifications, the CGA reduces to a truck-only version that is consistent with the problem introduced in \citet{khameneh2023non}.

\subsection{Numerical results}
This section presents the numerical results obtained under two operational settings: (i)~truck-only collection and (ii)~drone-aided collection. Evaluating both scenarios enables us to assess not only the effectiveness of the CGA when solving the problem considered by \citet{khameneh2023non}, but also the operational benefits achievable when drones are integrated into the blood collection system as proposed in this study. The results are organized in Tables~\ref{Table:NoDroneResults} and~\ref{Table:WithDroneResults}. The former reports the outcomes obtained when only trucks are available for blood collection, while the latter presents the results when each truck is equipped with a drone, enabling aerial pickups and synchronized ground–air operations. For all algorithms in both tables, the reported values correspond to the average performance over ten independent runs.

\input{Tables/Table_NoDroneResults}
\input{Tables/Table_WithDroneResults}

In Table~\ref{Table:NoDroneResults}, the first four columns describe the characteristics of each instance, including its class, identifier, number of DSs, and fleet size. The remaining columns are grouped by the solution approach. Within the \textit{truck-only model}~(TOM) section, the column \textit{Obj} reports the best objective value found by Gurobi when solving the mathematical model developed by \citet{khameneh2023non}, while \textit{UB} provides the best upper bound obtained during the search. The column \textit{Time~(sec)} denotes the total runtime in seconds, and \textit{Gap~(\%)} indicates the optimality gap computed as $100 \times (UB - Obj)/Obj$. For algorithms, the reported \textit{Obj} values correspond to the average objective function value over the ten runs, and \textit{Time (sec)} refers to the average computation time. Since Gurobi fails to obtain tight upper bounds in most cases, we use its best solution for performance comparison. Accordingly, the column \textit{Imp~(\%)} measures algorithm performance relative to the solution obtained by Gurobi, calculated as $100 \times (Z_A - Z_G)/Z_G$, where $Z_A$ and $Z_G$ denote the objective values obtained by the corresponding algorithm and Gurobi, respectively. When this metric takes negative values, it means that the solution found by Gurobi outperforms that of the corresponding algorithm. The \textit{CV} column presents the coefficient of variation of the objective values over the ten runs. Finally, the \textit{SR} column represents the success rate, which measures the proportion of runs for which a qualified solution is obtained. A qualified solution is defined as one whose objective value is at most 2\% worse than the best-known value achieved across all methods. 

Table~\ref{Table:WithDroneResults} is organized in the same manner as Table~\ref{Table:NoDroneResults}, yet reflects the outcomes obtained when aerial pickups are permitted. Two differences distinguish this table from the truck-only case. First, the section labeled \textit{Drone-aided model} presents the results obtained by solving the DAM with Gurobi rather than the TOM. Second, each algorithm section includes an additional metric, denoted $\textit{DG}~(\%)$, which measures the percentage gain achieved through drone integration.
This metric is computed as $ 100 \times (Z_A^{\text{best}} - Z_{\text{truck-only}}^{\text{best}})/ {Z_{\text{truck-only}}^{\text{best}}}$, where $Z_A^{\text{best}}$ denotes the best objective value obtained by the corresponding algorithm in the drone-aided setting, and $Z_{\text{truck-only}}^{\text{best}}$ represents the best objective value achieved across all methods when only trucks are available. 

The results indicate that Gurobi is able to solve both the TOM and the DAM to optimality for instances with up to 9~DSs (instances~1–4 in class~C1). In this range, CGA also finds the optimal solution in every run for both problem settings, while requiring substantially less computation time. Under the truck-only configuration, CGA finishes in only 1.80~seconds on average, compared with 3253.00~seconds for Gurobi, and this advantage persists in the drone-aided setting, with average runtimes of 2.02~seconds versus 2152.75~seconds. For these small instances, HGA performs similarly well, reaching the optimal solution in all runs; however, its computation times are higher than those of CGA, averaging 10.39 seconds under the truck-only configuration and 10.31 seconds in the drone-aided setting. IWO performs competitively under the truck-only configuration and, except for a single run of instance~4, it consistently reaches the optimal solution for the TOM. However, it becomes less stable in the drone-aided setting; while optimal solutions are still achieved for the DAM, they are obtained intermittently rather than in all ten runs. Despite these fluctuations, IWO remains significantly faster than Gurobi, with average runtimes of 3.50~seconds (truck-only) and 3.33~seconds (drone-aided).

For the larger instances in class C1 (instances 5–10), Gurobi can no longer provide competitive solutions, even after running for the full 12-hour computational budget. On average, its best truck-only and drone-aided solutions are 11.14\% and 9.28\% worse than those obtained by CGA. Compared with the two benchmark algorithms, CGA also delivers better and more reliable results. Under the truck-only configuration, its objective values exceed those of HGA and IWO by 2.17\% and 6.00\%, respectively, and it maintains the highest reliability, achieving qualified solutions in 92\% of runs. By contrast, HGA succeeds in only 53\% of the runs and IWO in 22\%. When drones are allowed, CGA maintains its advantage, outperforming HGA and IWO by 3.39\% and 7.97\%, and delivering qualified solutions in 90\% of runs, compared with only 17\% for HGA and 7\% for IWO. 

Moving to classes~C2 and C3, which contain medium to large-scale instances, the advantages of CGA become increasingly evident. When trucks handle all pickups, CGA not only outperforms Gurobi but also remains the most dependable solution approach. In class~C2, it improves Gurobi’s best incumbent by an average of 43.78\% and attains qualified solutions in 84\% of the runs, while HGA and IWO succeed in only 17\% and 1\%, respectively. This trend intensifies in class~C3, where CGA improves the solver’s best solutions by an average of 58.83\% and remains the most efficient approach, achieving a 87\% success rate compared with 15\% for HGA and 0\% for IWO. The performance gap becomes even more critical when drones are integrated. In class~C2, CGA surpasses Gurobi’s best incumbents by 41.83\% on average and remains stable with an 83\% success rate. By contrast, HGA is only able to produce qualified solutions for the smallest instance (with 23 DSs), and IWO fails to do so altogether. In class~C3, Gurobi fails to scale to the complexity of the problem in instances 6, 8, and 10, returning no feasible solution. Moreover, neither HGA nor IWO produces a single qualified solution across all runs. CGA, however, continues to perform consistently, achieving an 85\% success rate and improving Gurobi’s best incumbents by 147.36\% on average.

\begin{figure}[h]
    \centering
    \captionsetup{justification=centering} 
    \begin{subfigure}[b]{0.48\textwidth}
        \centering
        \includegraphics[scale=0.53]{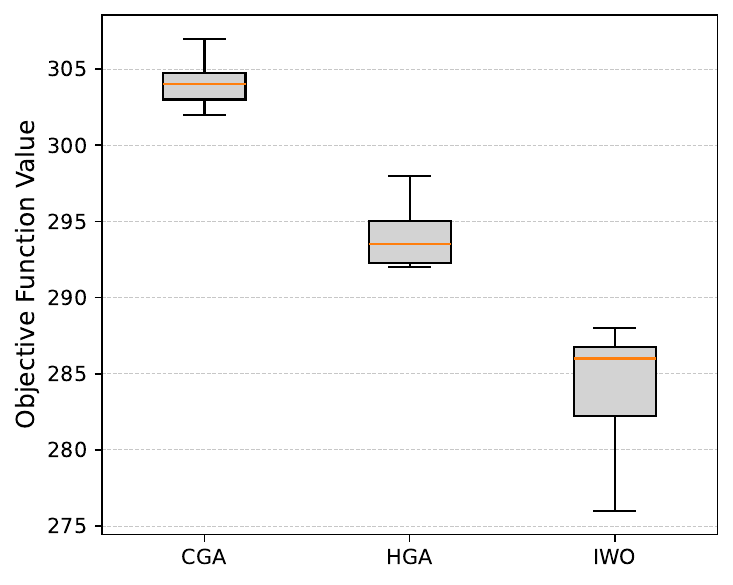}
        \caption{Truck-only system}
        \label{fig:BoxPlot_truck} 
    \end{subfigure}  
    \hfill 
    \begin{subfigure}[b]{0.48\textwidth}
        \centering
        \includegraphics[scale=0.53]{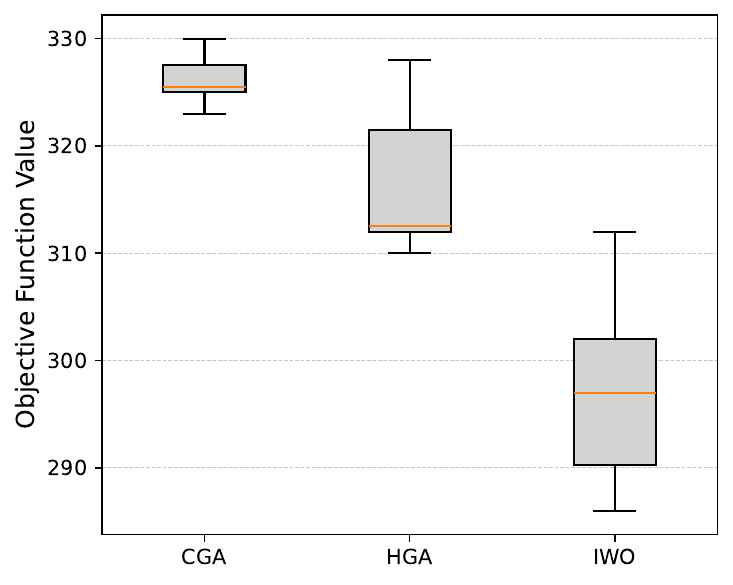}
        \caption{Drone-aided system}
        \label{fig:BoxPlot_drone} 
    \end{subfigure} 
    \vspace{0.25cm} 
    \caption{Box plots of the objective values for instance 1 of class C2.}
\end{figure}

To further illustrate the comparative performance of the algorithms, Figs.~\ref{fig:BoxPlot_truck}~and~\ref{fig:BoxPlot_drone} present box plots of the objective values obtained by CGA, HGA, and IWO for instance~1 of class~C2 under the truck-only and drone-aided settings, respectively. In the truck-only case, CGA exhibits very stable behavior, with a very narrow interquartile range and a CV of 0.47\%, compared with 0.72\% for HGA and 1.30\% for IWO. When drones are integrated, the performance gap widens and the behavior of the benchmark algorithms becomes less stable. Although HGA and IWO obtain higher objective values than in the truck-only case, their variability increases sharply, as reflected by larger interquartile ranges and CV exceeding 2\% and 3\%, respectively. In contrast, CGA maintains a tight distribution and low variability, with a CV of 0.62\%, highlighting its superior robustness.

\begin{figure}[h]
    \hspace{-0.5cm}
    \centering
    \captionsetup{justification=centering} 
    \begin{subfigure}[b]{0.50\textwidth}
        \centering
        \includegraphics[scale=0.43]{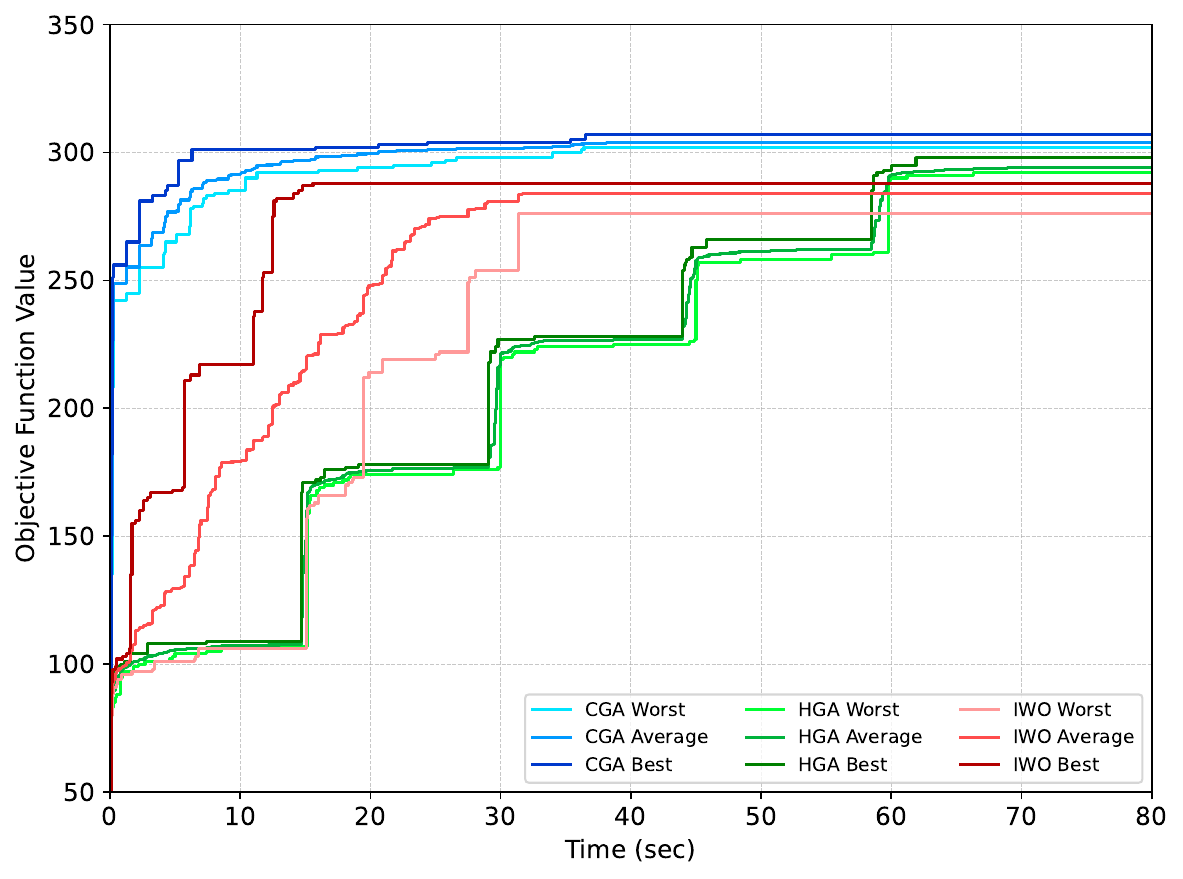}
        \caption{Truck-only system}
          \label{fig:ConvCurve_truck} 
    \end{subfigure}  
    \hfill 
    \begin{subfigure}[b]{0.50\textwidth}
        \centering
        \includegraphics[scale=0.43]{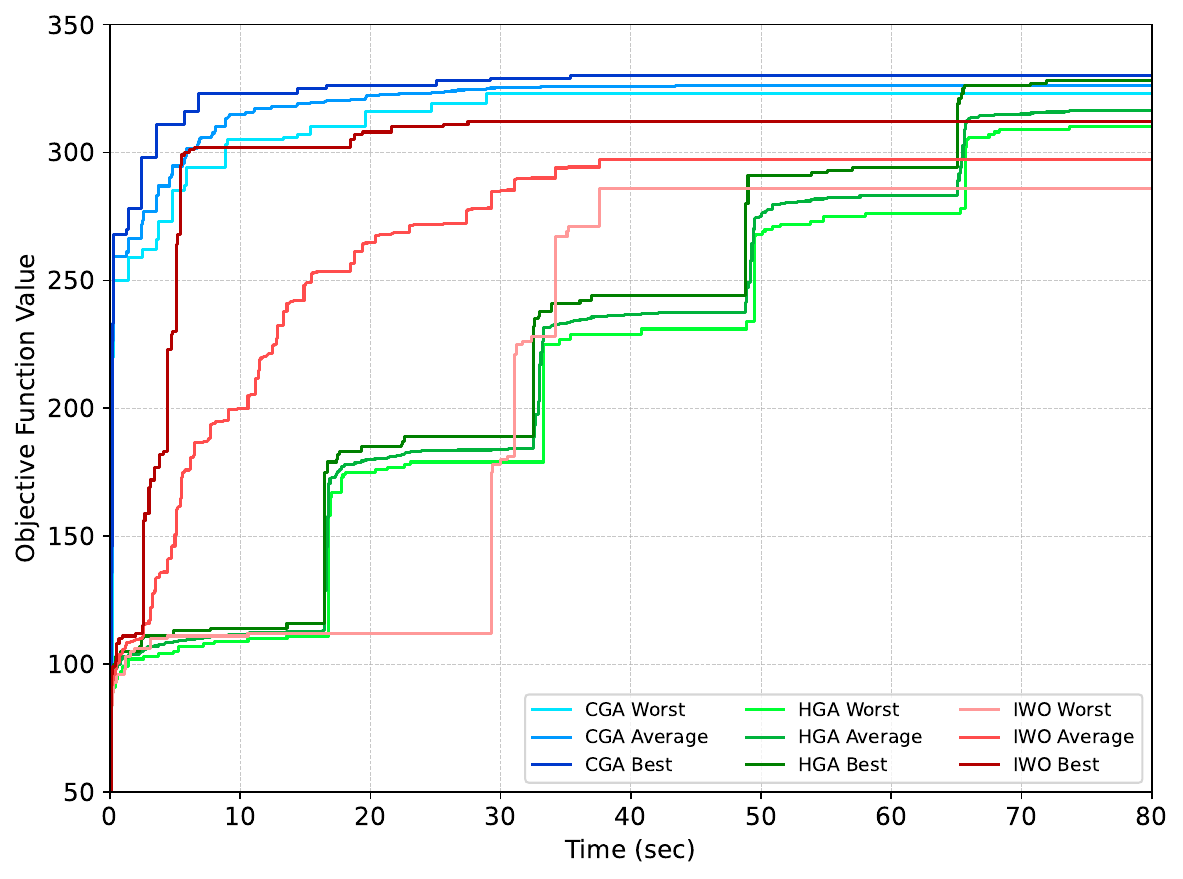}
        \caption{Drone-aided system}
        \label{fig:ConvCurve_drone} 
    \end{subfigure} 
    \vspace{0.25cm} 
    \caption{Convergence curves for instance~1 of class~C2.} 
\end{figure}

The convergence curves shown in Figs.~\ref{fig:ConvCurve_truck}~and~\ref{fig:ConvCurve_drone} further support these observations. These figures report the best, average, and worst objective values over time across ten independent runs for each algorithm under the truck-only and drone-aided settings. In both settings, CGA converges to higher-quality solutions compared to both benchmark algorithms, with a small gap between its best and worst trajectories, indicating strong robustness to stochastic effects. By contrast, HGA and IWO exhibit wider gaps between their best and worst curves, particularly in the drone-aided case. These results further underscore the value of CGA as a robust and reliable solution approach.

Finally, the results clearly demonstrate the benefits of integrating drones into the blood collection system. When evaluated using CGA, the best-performing solution approach, the inclusion of drones yields average DG values of 5.52\% in class~C1, 6.49\% in class~C2, and 4.81\% in class~C3 relative to the corresponding truck-only best solutions. These improvements translate into an average of 8.5, 30.5, and 37.5 additional blood units delivered to the processing center in classes~C1, C2, and~C3, respectively. Such gains are particularly important in blood supply operations, where even small increases in collected units can have meaningful implications for patient care and healthcare system responsiveness.

\subsection{Sensitivity analysis}
This section examines how variations in key operational factors influence the overall blood collection performance. The analyses are conducted on the instance exhibiting the highest DG in class~C1 (i.e., instance~8), and the CGA is used to obtain solutions in all experiments.

\subsubsection{Effect of fleet size on the collection system}

\begin{figure}[h]
    \hspace{-0.25cm}\makebox[\linewidth][c]{%
        \begin{subfigure}[t]{0.42\paperwidth}
            \centering
            \includegraphics[width=\linewidth]{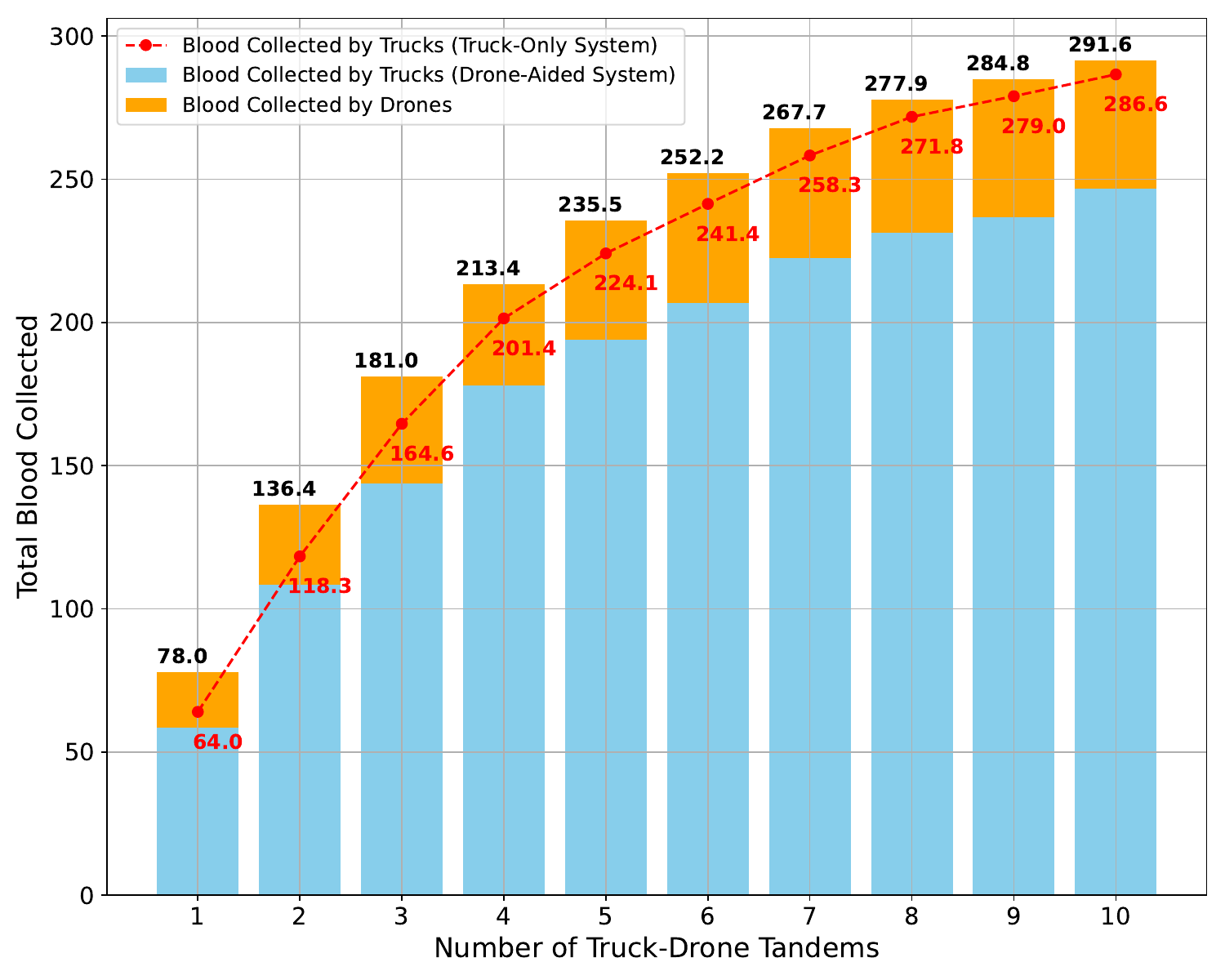}
            \caption{Total viable blood collected by each system.}
            \label{fig:platelet_production}
        \end{subfigure}%
         
         \hspace{0.30cm}
         
        \begin{subfigure}[t]{0.42\paperwidth}
            \centering
            \includegraphics[width=\linewidth]{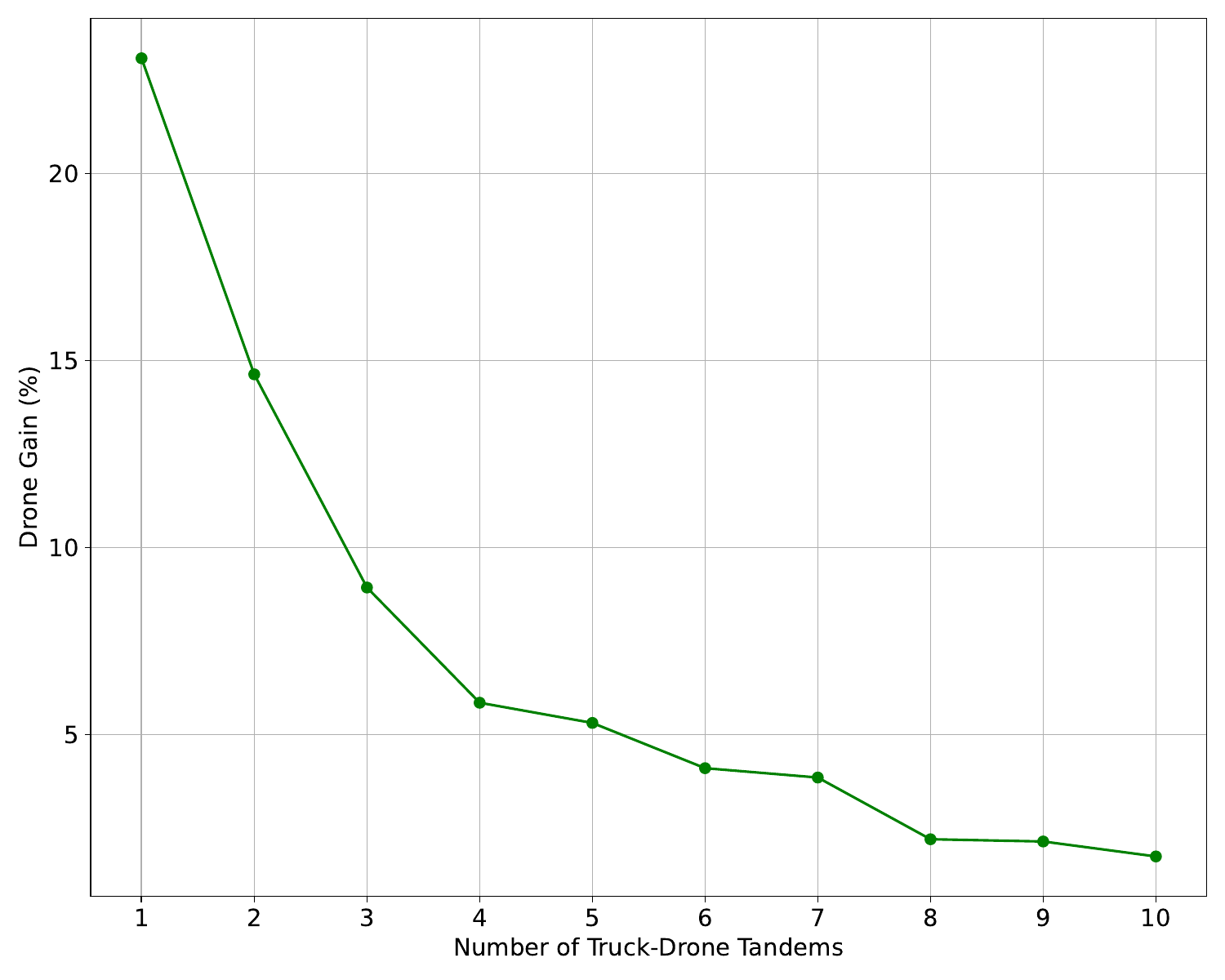}
            \caption{Performance improvement with drone integration.}
            \label{fig:percentage_improvement}
        \end{subfigure}%
    }
    \vspace{0.5cm}
    \caption{Comparison of drone-aided and truck-only systems.}
    \label{fig:Sensitivity_PartOne}
\end{figure}

Fig.~\ref{fig:Sensitivity_PartOne} compares the drone-aided and truck-only systems under varying numbers of truck–drone tandems. Fig.~\ref{fig:platelet_production} reports the total viable blood collected in each system, averaged over ten CGA runs, while Fig.~\ref{fig:percentage_improvement} quantifies the relative improvement achieved through drone integration using the DG metric defined in the previous section.

The results show that the drone-aided system consistently outperforms the truck-only system across all fleet sizes. The relative benefit of drone integration is most pronounced in resource-limited scenarios, with DG exceeding 23\% when only one tandem is available. As the number of tandems increases, DG decreases steadily, dropping below 9\% at three tandems and falling to less than 2\% when ten tandems are deployed. This trend is expected, as increasing the fleet size allows trucks to cover more sites and collect larger volumes of blood, thereby reducing the marginal contribution of drones. Nevertheless, the drone-aided system remains superior throughout, underscoring the operational value of drone integration in increasing system efficiency under varying resource levels.

\subsubsection{Effect of drone-related parameters on the collection system}
To evaluate the impact of drone-related parameters on the performance of the system, we conducted a sensitivity analysis on the payload capacity, battery endurance, and speed of drones. In each experiment, one parameter is varied across three levels—low, medium, and high—while the remaining parameters are held at their medium levels. The specific values used for each level are summarized in Table~\ref{Table_SensitivityLevels}.

\input{Tables/Table_SensAnalLevels}

The results are illustrated in Fig.~\ref{fig:Sensitivity_DroneOnlyHeatMaps}, where each heat map shows the blood collected by drones for various levels of the parameter tested, across tandem counts ranging from one to four. An analysis of the solutions reveals the following key findings:

\begin{figure}[h]
   \hspace*{-0.7cm}
   \includegraphics[scale=0.45]{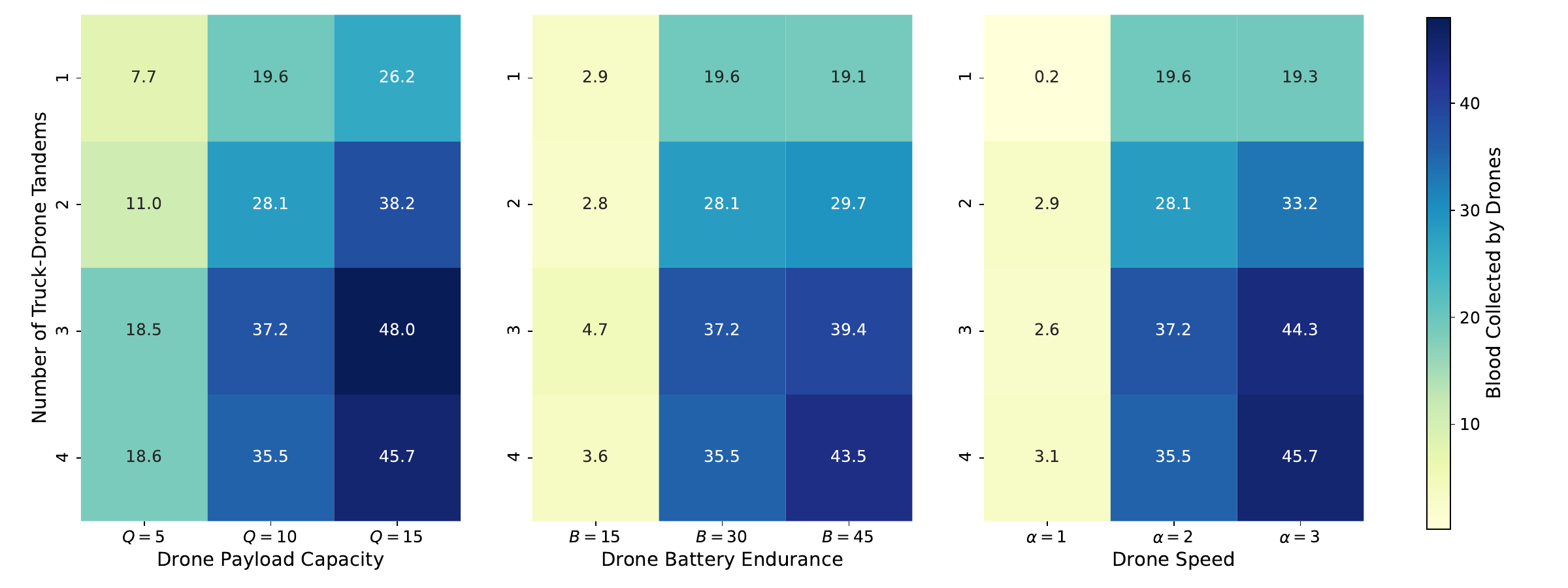}
   \vspace{0.02cm}
   \caption{Impact of drone-related parameters on the number of blood units collected by drones.}
   \label{fig:Sensitivity_DroneOnlyHeatMaps}
\end{figure}

\begin{itemize}[label={$\bullet$}]

    \item In general, equipping the fleet with more advanced drones in terms of payload capacity, battery endurance, or speed leads to a clear increase in the amount of blood collected by drones. However, the magnitude of these improvements varies depending on the parameter considered and the number of available truck–drone tandems.

    \item Speed has the most restrictive effect. When drones operate at the same speed as trucks ($\alpha=1$), their contribution is negligible, remaining below 3.1 units across all fleet sizes. Increasing the speed ratio to $\alpha=2$ leads to a sharp increase in drone effectiveness, with collected volumes reaching values between 19.6 and 37.2 units depending on the fleet size. Further increasing the speed to $\alpha=3$ yields additional gains, particularly for larger fleets, where drone-collected blood reaches up to 45.7 units.
    
    \item Battery endurance is the second most impactful parameter. With a limited endurance of $B=15$, drones collect fewer than 4.7 units regardless of fleet size, indicating severe operational constraints. Increasing $B$ from 15 to
    30 leads to substantial improvements in drone performance. For instance, at three tandems, the amount of blood collected by drones increases from 4.7 to 37.2 units. Although the gain from increasing $B$ from 30 to 45 is less pronounced, it remains meaningful; at four tandems, the drone-collected volume increases by an additional 8 units.

    \item Payload capacity has the least restrictive effect when set to its lowest level, compared to battery endurance and speed. Indeed, when $B$ and $\alpha$ are at their minimum values, drones are almost ineffective, regardless of their payload capacity. However, once drones have adequate battery life and speed, even a small payload capacity ($Q=5$) allows them to contribute meaningfully.  Increasing $Q$ to its highest level ($Q=15$) leads to the largest amount of blood collected by drones, exceeding the contributions achieved when either $B$ or $\alpha$ is set to its maximum level.
\end{itemize}

\section{Conclusion} \label{Sec7}
This study addressed the drone-aided blood collection routing problem, which extends the truck-only blood collection routing problem by integrating synchronized truck–drone operations. This collaborative framework increases the flexibility of the collection fleet by allowing trucks and drones to perform pickups in parallel, thereby increasing the number of blood units that reach the processing center while they remain viable for platelet extraction. To the best of our knowledge, this study is the first to consider coordinated truck–drone operations in the blood collection routing problem.

We formulated the problem as an MILP model that maximizes the total number of viable blood units collected by the fleet by jointly optimizing truck routes, drone sorties, pickup schedules, and timing decisions. Due to the model’s computational intractability for Gurobi, we proposed a column generation algorithm that decomposes the problem into a master problem selecting truck–drone tours and a pricing subproblem generating promising ones. We developed an efficient memetic algorithm for addressing the pricing subproblem that generates high-quality tours within reasonable computation times. We evaluated our proposed column generation algorithm against two state-of-the-art metaheuristics proposed in the literature, namely the hybrid genetic algorithm and invasive weed optimization. The numerical results demonstrated that our proposed algorithm consistently outperforms the benchmark methods in both the truck-only and drone-aided settings, while remaining robust and computationally efficient across a wide range of instance sizes. Beyond algorithmic performance, the results clearly demonstrate the operational value of drone integration into the blood collection fleet. Across all tested scenarios, the drone-aided system consistently collected more viable blood units than the truck-only system. These improvements are particularly important in healthcare systems, where every additional unit of collected blood directly supports life-saving treatments. Finally, we conducted a sensitivity analysis to examine how variations in key operational factors, such as fleet size and drone capabilities, influence overall performance. These analyses revealed that drone integration is most impactful in resource-limited scenarios and that faster, longer-range, and higher-capacity drones provide substantially greater operational benefits.

Future work could broaden this research in several directions. One promising extension is to incorporate uncertainty in travel times, donor arrivals, and donation volumes, leading to stochastic or robust variants of the problem that better reflect real-world operating conditions. Another direction is to incorporate speed- and payload-dependent energy consumption for drones, enabling more realistic routing and scheduling decisions. In addition, extending the framework to multi-processing-center networks and to multiple blood products with different perishability profiles would further enhance its applicability in large-scale blood supply chains. Finally, integrating strategic decisions such as drone base locations, fleet sizing, and investment planning with the operational routing model would provide valuable insights for policymakers and healthcare organizations seeking to deploy drone-supported blood collection systems.

\appendix
\renewcommand{\thesection}{Appendix \Alph{section}}
\section{Road-travel time matrix} \label{sec:appTimeMatrix}

Table~\ref{tab:appendix_truck_times} reports the road-travel time matrix for the instance presented in Section~\ref{Sec1}, where each entry represents the time (in minutes) required for a truck to travel between two locations. For example, the truck travel time from DS~1 to DS~3 is 110 minutes.

\begin{table}[h]
\centering
\caption{Road-travel times between locations (minutes)}
\label{tab:appendix_truck_times}
\begin{tblr}{
  cells = {c},
  cell{1}{2} = {c=5}{},
  hline{1,3,8} = {-}{},
  hline{2} = {2-6}{},
}
     & To   &      &      &      &     \\
From & DS 1 & DS 2 & DS 3 & DS 4 & PC  \\
DS 1 & 0    & 30   & 110  & 30   & 130 \\
DS 2 & 30   & 0    & 120  & 30   & 120 \\
DS 3 & 110  & 120  & 0    & 90   & 140 \\
DS 4 & 30   & 30   & 90   & 0    & 150 \\
PC   & 130  & 120  & 140  & 150  & 0   
\end{tblr}
\end{table}

\section{Local-search operators} \label{sec:appLSOperators}
The sixteen local-search operators used in the education phase are as follows.  Note that, in all operators, the involved DSs are selected randomly from the current route unless otherwise stated.

{\small
\begin{itemize}

\item[--] $M_1$: Relocate a truck-only DS~$i$ (a DS without any drone launch or retrieval activity) after a DS~$j$ in the truck route.

\item[--] $M_2$: Relocate two consecutive truck-only DSs $i_1$ and $i_2$ after a DS $j$ in the truck route, preserving the order $i_1$ followed by $i_2$.

\item[--] $M_3$: Relocate two consecutive truck-only DSs $i_1$ and $i_2$ after a DS $j$ in the truck route, reversing the order to $i_2$ followed by $i_1$.

\item[--] $M_4$: Swap two DSs $i$ and $j$ in the truck route.

\item[--] $M_5$: Swap two consecutive DSs $i_1$ and $i_2$ in the truck route with another DS $j$, where $i_2$ has no drone launch or retrieval activity.

\item[--] $M_6$: Swap two consecutive DSs $i_1$ and $i_2$ in the truck route with two other consecutive DSs $j_1$ and $j_2$.

\item[--] $M_7$: Select two pairs of consecutive DSs $(i_1,i_2)$ and $(j_1,j_2)$ in the truck route, and reconnect them as $(i_1,j_1)$ and $(i_2,j_2)$.

\item[--] $M_8$: Select two pairs of consecutive DSs $(i_1,i_2)$ and $(j_1,j_2)$ in the truck route, and reconnect them as $(i_1,j_2)$ and $(j_1,i_2)$.

\item[--] $M_9$: Swap a drone-served DS $i$ with a truck-served DS $j$, where $j$ is neither the launch nor the rendezvous DS of the drone visiting $i$, nor any DS between them.

\item[--] $M_{10}$: In a drone sortie $(i,n,j)$, swap the launch DS $i$ with the drone-served DS $n$.

\item[--] $M_{11}$: In a drone sortie $(i,n,j)$, swap the drone-served DS $n$ with the rendezvous DS $j$.

\item[--] $M_{12}$: In a drone sortie $(i,n,j)$, swap the launch DS $i$ with the rendezvous DS $j$.

\item[--] $M_{13}$: Insert a new drone sortie $(i, n, j)$ by selecting a DS $n$ that is neither a launch nor a rendezvous node, and two truck-served DSs $i$ and $j$ ($i$ preceding $j$), provided that no other drone launch or retrieval occurs between $i$ and $j$.

\item[--] $M_{14}$: Remove a drone-served DS $n$ from its drone sortie and insert it as a truck-served DS between two consecutive DSs in the truck route.

\item[--] $M_{15}$:
Select two drone sorties $(i_1,n_1,j_1)$ and $(i_2,n_2,j_2)$ and swap the drone-served DSs $n_1$ and $n_2$, resulting in two new sorties $(i_1,n_2,j_1)$ and $(i_2,n_1,j_2)$.

\item[--] $M_{16}$: Reassign a drone-served DS $n$ to a new drone sortie by selecting new launch and rendezvous DSs.

\end{itemize}
}

\bibliographystyle{unsrtnat}
\bibliography{references} 

\end{document}

%% file: Tables/Table_Intro.tex
\begin{table}[h]
\centering
\caption{ Donation counts and completion times at each donation site.}
\label{Table:DonationPlan}
\scalebox{1.0}{
\def\arraystretch{1}
\begin{tblr}{
  colspec = {
    Q[l,wd=2.5cm]  
    Q[c,wd=1.5cm]  
    Q[c,wd=1.5cm]
    Q[c,wd=1.5cm]
    Q[c,wd=1.5cm]
  },
  cell{1}{1} = {r=2}{},
  cell{1}{2} = {c=4}{},
  hline{1,3,9} = {-}{},
  hline{2} = {2-5}{},
}
{Donation\\completion time} & Number of donations &            &             &            \\
              & DS 1                & DS 2       & DS 3        & DS 4      \\
10:30         & 0                   & 0          & 3           & 4          \\
11:00         & 5                   & 6          & 0           & 2          \\
11:30         & 0                   & 0          & 4           & 0          \\
12:00         & 0                   & 0          & 1           & 0          \\
12:30         & 0                   & 0          & 5           & 0          \\
\textbf{Total}  & \textbf{5}          & \textbf{6} & \textbf{13} & \textbf{6} 
\end{tblr}}
\end{table}

%% file: Tables/Table_Notation.tex
\begin{table}[h]
\caption{Notation}
\label{Table_Notation}
\centering
\renewcommand{\arraystretch}{1.35} 
\begin{tabular}{l@{\hspace{0.10cm}}p{13cm}}
\hline
\textbf{Sets} & \\
$\mathcal{V}=\{1,2,\dots,N\}$ & Set of DSs - indexed by $i,j,n,n^\prime$ \\
$\mathcal{V}_0=\mathcal{N} \cup \{0\}$ & Set of DSs and PC - indexed by $i,j,n,n^\prime$ \\
$E$ & Set of edges connecting the nodes in $\mathcal{V}_0$ - indexed by $(i,j)$ \\
$F$ & Set of truck-drone tandems - indexed by $k$ \\
$\mathcal{T}=\{1,2,\dots,L\}$ & Set of routes performed by a truck-drone tandem in a single day - indexed by $m$ \\
$T_i=\{a_i,\dots,b_i\}$ & Set of donation completion times at DS~$i$ - indexed by $t$ \\
$\mathcal{D}$ & Set of feasible sorties - indexed by $(i,n,j)$\\
\vspace{0.15pt} \\

\textbf{Parameters} & \\
$t^v_{ij}$ & Travel time on edge~$(i,j)$ by a truck \\
$t^d_{ij}$ & Travel time on edge~$(i,j)$ by a drone  \\
$Q$ & Payload capacity of each drone \\
$B$ & Battery endurance of each drone \\
$a_i$ & Earliest blood donation completion time at DS~$i$ \\
$b_i$ & Latest blood donation completion time at DS~$i$ \\
$T_s$ & Start time of planning horizon \\
$T_e$ & End time of planning horizon \\
$d_{it}$ & Number of blood donations completed at DS~$i$ at time~$t$  \\
$K$ & Processing time limit \\
$H$ & A large enough number \\

\vspace{0.15pt} \\

\textbf{Decision variables} & \\
$x_{itmk}$ & 1 if blood donations completed at DS~$i$ at time~$t$ are picked up in route $m$ of truck~$k$; 0 otherwise \\
$f_{itmk}$ & 1 if blood donations completed at DS~$i$ at time~$t$ are picked up in route $m$ of drone~$k$; 0 otherwise \\
$y_{ijmk}$ &  1 if node $j$ is visited immediately after node $i$ in route $m$ of truck $k$; 0 otherwise \\
$h_{injmk}$ & 1 if sortie $(i,n,j)$ is performed in route $m$ of drone $k$; 0 otherwise \\
$r_{imk}$ & 1 if a sortie is initiated from node $i$ in route $m$ of truck-drone tandem $k$; 0 otherwise. \\ 
$z_{mk}$ & 1 if route $m$ of truck-drone tandem $k$ is performed; 0 otherwise \\
$v_{imk}$ & Departure time from DS~$i$ in route $m$ of truck-drone tandem $k$ \\
$e_{mk}$ & Completion time of route $m$ of truck-drone tandem $k$ \\
\hline
\end{tabular}
\end{table}

%% file: Tables/Table_Taguchi.tex
\begin{table}[h]
\caption{Tuned parameters.}
\label{tab:Taguchi}
\centering
\small
\begin{tblr}{
  colspec = {p{2cm} p{1.5cm} p{1.5cm} p{1.5cm}},
  hline{1-2,8} = {-}{},
}
Parameter          & Level 1 & Level 2 & Level 3 & Optimal level   \\
\texttt{poolSize}  & $3\times 10^{5}$  & $5\times 10^{5}$  & $7\times 10^{5}$  & $5\times 10^{5}$        \\
\texttt{genLimit}  & 600               & 800               & 1000              & 1000                    \\
\texttt{popSize}   & 100               & 120               & 140               & 140                     \\
\texttt{refIter}   & 2                 & 8                 & 16                & 16                      \\
\texttt{muRate}    & 0.1               & 0.3               & 0.6               & 0.3                     \\
\texttt{muVolume}  & 0.1               & 0.2               & 0.3               & 0.3            
\end{tblr}
\end{table}

%% file: Tables/Table_NoDroneResults.tex
\begin{table}[p]
\centering
\def\arraystretch{1.5}

\rotatebox{90}{%
\begin{minipage}{\textheight} 

\caption{Computational results for the truck-only blood collection system.}
\label{Table:NoDroneResults}

\scalebox{0.69}{%
\sffamily

\hspace{-0cm}\begin{tabular}{lllllllllllllllllllllllllll} 
\hline
\multirow{2}{*}{Class} & \multirow{2}{*}{Inst.} & \multirow{2}{*}{DSs} & \multirow{2}{*}{$|F|$} &  & \multicolumn{4}{l}{Truck-only model (MILP)}                                                                             &  & \multicolumn{5}{l}{CGA proposed in our study}                                                                                                                                     &  & \multicolumn{5}{l}{HGA proposed by \cite{khameneh2023non}}                                                                                                                                     &  & \multicolumn{5}{l}{IWO proposed by \cite{khameneh2023non}}                                                                                                                                      \\ 
\cline{6-9}\cline{11-15}\cline{17-21}\cline{23-27}
                       &                        &                       &                                        &  & Obj             & UB               & \begin{tabular}[c]{@{}l@{}}Time\\(sec)\end{tabular} & \begin{tabular}[c]{@{}l@{}}Gap\\(\%)\end{tabular} &  & Obj             & \begin{tabular}[c]{@{}l@{}}Time\\(sec)\end{tabular} & \begin{tabular}[c]{@{}l@{}}Imp\\(\%)\end{tabular} & CV              & SR            &  & Obj             & \begin{tabular}[c]{@{}l@{}}Time\\(sec)\end{tabular} & \begin{tabular}[c]{@{}l@{}}Imp\\(\%)\end{tabular} & CV              & SR            &  & Obj             & \begin{tabular}[c]{@{}l@{}}Time\\(sec)\end{tabular} & \begin{tabular}[c]{@{}l@{}}Imp\\(\%)\end{tabular} & CV              & SR             \\ 
\hline
C1                     & 1                      & 3                     & 1                                      &  & 28              & 28               & 11                                                  & 0.00                                              &  & 28.0            & 0.55                                                & 0.00                                              & 0.0000          & 1.0           &  & 28.0            & 5.11                                                & 0.00                                              & 0.0000          & 1.0           &  & 28.0            & 0.06                                                & 0.00                                              & 0.0000          & 1.0            \\
                       & 2                      & 5                     & 1                                      &  & 37              & 37               & 85                                                  & 0.00                                              &  & 37.0            & 1.02                                                & 0.00                                              & 0.0000          & 1.0           &  & 37.0            & 6.38                                                & 0.00                                              & 0.0000          & 1.0           &  & 37.0            & 1.90                                                & 0.00                                              & 0.0000          & 1.0            \\
                       & 3                      & 7                     & 2                                      &  & 43              & 43               & 9450                                                & 0.00                                              &  & 43.0            & 2.09                                                & 0.00                                              & 0.0000          & 1.0           &  & 43.0            & 14.44                                               & 0.00                                              & 0.0000          & 1.0           &  & 43.0            & 4.87                                                & 0.00                                              & 0.0000          & 1.0            \\
                       & 4                      & 9                     & 2                                      &  & 60              & 60               & 3466                                                & 0.00                                              &  & 60.0            & 3.55                                                & 0.00                                              & 0.0000          & 1.0           &  & 60.0            & 15.61                                               & 0.00                                              & 0.0000          & 1.0           &  & 59.7            & 7.16                                                & -0.50                                             & 0.0081          & 1.0            \\
                       & 5                      & 11                    & 2                                      &  & 124             & 198              & 43200                                               & 59.68                                             &  & 125.0           & 14.60                                               & 0.81                                              & 0.0000          & 1.0           &  & 124.9           & 20.20                                               & 0.73                                              & 0.0025          & 1.0           &  & 122.7           & 8.19                                                & -1.05                                             & 0.0154          & 0.9            \\
                       & 6                      & 13                    & 3                                      &  & 101             & 250              & 43200                                               & 147.53                                            &  & 105.6           & 6.45                                                & 4.55                                              & 0.0066          & 1.0           &  & 105.6           & 33.56                                               & 4.55                                              & 0.0049          & 1.0           &  & 101.6           & 7.50                                                & 0.59                                              & 0.0124          & 0.3            \\
                       & 7                      & 15                    & 3                                      &  & 166             & 366              & 43200                                               & 120.48                                            &  & 181.2           & 30.94                                               & 9.16                                              & 0.0081          & 0.9           &  & 175.1           & 40.79                                               & 5.48                                              & 0.0032          & 0.0           &  & 169.6           & 12.24                                               & 2.17                                              & 0.0176          & 0.0            \\
                       & 8                      & 17                    & 3                                      &  & 157             & 359              & 43200                                               & 128.66                                            &  & 164.6           & 19.65                                               & 4.84                                              & 0.0152          & 0.6           &  & 161.8           & 44.73                                               & 3.06                                              & 0.0136          & 0.3           &  & 158.3           & 14.88                                               & 0.83                                              & 0.0219          & 0.1            \\
                       & 9                      & 19                    & 4                                      &  & 219             & 507              & 43200                                               & 131.51                                            &  & 272.8           & 23.61                                               & 24.57                                             & 0.0034          & 1.0           &  & 254.2           & 61.30                                               & 16.07                                             & 0.0052          & 0.0           &  & 247.6           & 12.95                                               & 13.06                                             & 0.0257          & 0.0            \\
                       & 10                     & 21                    & 4                                      &  & 256             & 571              & 43200                                               & 123.05                                            &  & 314.6           & 22.09                                               & 22.89                                             & 0.0075          & 1.0           &  & 313.4           & 58.32                                               & 22.42                                             & 0.0074          & 0.9           &  & 288.1           & 26.74                                               & 12.54                                             & 0.0294          & 0.0            \\ 
\hline
\textbf{Average}       &                        &                       &                                        &  & \textbf{119.10} & \textbf{241.90}  & \textbf{27221.20}                                   & \textbf{71.09}                                    &  & \textbf{133.18} & \textbf{12.46}                                      & \textbf{6.68}                                     & \textbf{0.0041} & \textbf{0.95} &  & \textbf{130.30} & \textbf{30.04}                                      & \textbf{5.23}                                     & \textbf{0.0037} & \textbf{0.72} &  & \textbf{125.56} & \textbf{9.65}                                       & \textbf{2.76}                                     & \textbf{0.0131} & \textbf{0.53}  \\ 
\hline
C2                     & 1                      & 23                    & 5                                      &  & 211             & 530              & 43200                                               & 151.19                                            &  & 304.0           & 40.63                                               & 44.08                                             & 0.0047          & 1.0           &  & 294.1           & 73.57                                               & 39.38                                             & 0.0072          & 0.0           &  & 284.2           & 26.32                                               & 34.69                                             & 0.0130          & 0.0            \\
                       & 2                      & 25                    & 5                                      &  & 290             & 649              & 43200                                               & 123.79                                            &  & 356.9           & 47.74                                               & 23.07                                             & 0.0082          & 0.9           &  & 350.1           & 83.81                                               & 20.72                                             & 0.0107          & 0.2           &  & 339.2           & 43.54                                               & 16.97                                             & 0.0255          & 0.1            \\
                       & 3                      & 27                    & 6                                      &  & 382             & 740              & 43200                                               & 93.72                                             &  & 486.2           & 63.81                                               & 27.28                                             & 0.0092          & 0.7           &  & 456.9           & 114.86                                              & 19.61                                             & 0.0185          & 0.0           &  & 432.3           & 32.23                                               & 13.17                                             & 0.0353          & 0.0            \\
                       & 4                      & 29                    & 6                                      &  & 317             & 795              & 43200                                               & 150.79                                            &  & 445.3           & 70.57                                               & 40.47                                             & 0.0112          & 0.8           &  & 427.7           & 134.92                                              & 34.92                                             & 0.0145          & 0.0           &  & 405.2           & 56.59                                               & 27.82                                             & 0.0190          & 0.0            \\
                       & 5                      & 31                    & 6                                      &  & 274             & 744              & 43200                                               & 171.53                                            &  & 407.5           & 66.78                                               & 48.72                                             & 0.0090          & 0.8           &  & 405.8           & 131.79                                              & 48.10                                             & 0.0060          & 1.0           &  & 382.1           & 43.63                                               & 39.45                                             & 0.0101          & 0.0            \\
                       & 6                      & 33                    & 7                                      &  & 350             & 743              & 43200                                               & 112.29                                            &  & 475.8           & 49.71                                               & 35.94                                             & 0.0120          & 0.7           &  & 453.3           & 151.60                                              & 29.51                                             & 0.0152          & 0.0           &  & 429.1           & 43.13                                               & 22.60                                             & 0.0176          & 0.0            \\
                       & 7                      & 35                    & 7                                      &  & 338             & 1000             & 43200                                               & 195.86                                            &  & 538.5           & 79.65                                               & 59.32                                             & 0.0072          & 0.9           &  & 531.7           & 165.34                                              & 57.31                                             & 0.0068          & 0.3           &  & 492.1           & 47.61                                               & 45.59                                             & 0.0229          & 0.0            \\
                       & 8                      & 37                    & 7                                      &  & 363             & 855              & 43200                                               & 135.54                                            &  & 595.9           & 150.70                                              & 64.16                                             & 0.0109          & 0.7           &  & 569.0           & 185.66                                              & 56.75                                             & 0.0045          & 0.0           &  & 545.4           & 68.84                                               & 50.25                                             & 0.0202          & 0.0            \\
                       & 9                      & 39                    & 8                                      &  & 367             & 821              & 43200                                               & 123.71                                            &  & 512.5           & 144.94                                              & 39.65                                             & 0.0055          & 1.0           &  & 500.2           & 245.72                                              & 36.29                                             & 0.0109          & 0.2           &  & 442.5           & 69.81                                               & 20.57                                             & 0.0324          & 0.0            \\
                       & 10                     & 41                    & 8                                      &  & 429             & 1062             & 43200                                               & 147.55                                            &  & 665.5           & 123.50                                              & 55.13                                             & 0.0073          & 0.9           &  & 629.7           & 200.95                                              & 46.78                                             & 0.0130          & 0.0           &  & 602.4           & 93.87                                               & 40.42                                             & 0.0213          & 0.0            \\ 
\hline
\textbf{Average}       &                        &                       &                                        &  & \textbf{332.10} & \textbf{793.90}  & \textbf{43200.00}                                   & \textbf{140.60}                                   &  & \textbf{478.81} & \textbf{83.80}                                      & \textbf{43.78}                                    & \textbf{0.0085} & \textbf{0.84} &  & \textbf{461.85} & \textbf{148.82}                                     & \textbf{38.94}                                    & \textbf{0.0107} & \textbf{0.17} &  & \textbf{435.45} & \textbf{52.56}                                      & \textbf{31.15}                                    & \textbf{0.0217} & \textbf{0.01}  \\ 
\hline
C3                     & 1                      & 43                    & 9                                      &  & 370             & 1038             & 43200                                               & 180.54                                            &  & 580.9           & 102.58                                              & 57.00                                             & 0.0060          & 1.0           &  & 557.2           & 253.29                                              & 50.59                                             & 0.0106          & 0.0           &  & 513.7           & 82.67                                               & 38.84                                             & 0.0167          & 0.0            \\
                       & 2                      & 45                    & 9                                      &  & 342             & 1042             & 43200                                               & 204.68                                            &  & 636.9           & 200.65                                              & 86.23                                             & 0.0072          & 1.0           &  & 612.2           & 280.81                                              & 79.01                                             & 0.0121          & 0.0           &  & 584.6           & 101.77                                              & 70.94                                             & 0.0186          & 0.0            \\
                       & 3                      & 47                    & 10                                     &  & 321             & 980              & 43200                                               & 205.30                                            &  & 634.6           & 191.86                                              & 97.69                                             & 0.0077          & 0.9           &  & 620.5           & 319.66                                              & 93.30                                             & 0.0100          & 0.1           &  & 591.9           & 117.70                                              & 84.39                                             & 0.0146          & 0.0            \\
                       & 4                      & 49                    & 10                                     &  & 585             & 1248             & 43200                                               & 113.33                                            &  & 877.0           & 786.93                                              & 49.91                                             & 0.0088          & 0.8           &  & 855.1           & 345.79                                              & 46.17                                             & 0.0088          & 0.0           &  & 812.7           & 117.00                                              & 38.92                                             & 0.0170          & 0.0            \\
                       & 5                      & 51                    & 10                                     &  & 597             & 1314             & 43200                                               & 120.10                                            &  & 843.8           & 691.36                                              & 41.34                                             & 0.0086          & 0.7           &  & 832.8           & 346.82                                              & 39.50                                             & 0.0108          & 0.3           &  & 774.3           & 120.49                                              & 29.70                                             & 0.0236          & 0.0            \\
                       & 6                      & 53                    & 11                                     &  & 566             & 1356             & 43200                                               & 139.58                                            &  & 796.5           & 633.54                                              & 40.72                                             & 0.0085          & 0.9           &  & 776.0           & 415.59                                              & 37.10                                             & 0.0117          & 0.1           &  & 732.3           & 151.10                                              & 29.38                                             & 0.0213          & 0.0            \\
                       & 7                      & 55                    & 11                                     &  & 683             & 1464             & 43200                                               & 114.35                                            &  & 931.3           & 741.79                                              & 36.35                                             & 0.0060          & 1.0           &  & 921.4           & 403.24                                              & 34.90                                             & 0.0110          & 0.3           &  & 856.3           & 140.48                                              & 25.37                                             & 0.0141          & 0.0            \\
                       & 8                      & 57                    & 11                                     &  & 467             & 1255             & 43200                                               & 168.74                                            &  & 799.4           & 519.80                                              & 71.18                                             & 0.0078          & 0.9           &  & 789.0           & 418.32                                              & 68.95                                             & 0.0083          & 0.3           &  & 733.2           & 164.32                                              & 57.00                                             & 0.0153          & 0.0            \\
                       & 9                      & 59                    & 12                                     &  & 629             & 1481             & 43200                                               & 135.45                                            &  & 909.6           & 672.49                                              & 44.61                                             & 0.0084          & 0.7           &  & 897.3           & 469.65                                              & 42.66                                             & 0.0068          & 0.2           &  & 836.8           & 141.56                                              & 33.04                                             & 0.0132          & 0.0            \\
                       & 10                     & 61                    & 12                                     &  & 606             & 1630             & 43200                                               & 168.98                                            &  & 989.6           & 425.34                                              & 63.30                                             & 0.0092          & 0.8           &  & 972.8           & 496.55                                              & 60.53                                             & 0.0086          & 0.2           &  & 900.2           & 206.49                                              & 48.55                                             & 0.0121          & 0.0            \\ 
\hline
\textbf{Average}       &                        &                       &                                        &  & \textbf{516.60} & \textbf{1280.80} & \textbf{43200.00}                                   & \textbf{155.11}                                   &  & \textbf{799.96} & \textbf{496.63}                                     & \textbf{58.83}                                    & \textbf{0.0078} & \textbf{0.87} &  & \textbf{783.43} & \textbf{374.97}                                     & \textbf{55.27}                                    & \textbf{0.0099} & \textbf{0.15} &  & \textbf{733.60} & \textbf{134.36}                                     & \textbf{45.61}                                    & \textbf{0.0167} & \textbf{0.00}  \\
\hline
\end{tabular}}

\end{minipage}
}
\end{table}

%% file: Tables/Table_WithDroneResults.tex
\begin{table}[p]
\centering
\def\arraystretch{1.7}

\rotatebox{90}{%
\begin{minipage}{\textheight} 
\centering

\caption{Computational results for the drone-aided blood collection system.}
\label{Table:WithDroneResults}

\scalebox{0.61}{%

\sffamily
\hspace{-0cm}\begin{tabular}{llllllllllllllllllllllllllllll} 
\hline
\multirow{2}{*}{Class} & \multirow{2}{*}{Inst.} & \multirow{2}{*}{DSs} & \multirow{2}{*}{$|F|$} &  & \multicolumn{4}{l}{Drone-aided model (MILP)}                                                                                                                       &  & \multicolumn{6}{l}{CGA proposed in our study}                                                                                                                                                                                           &  & \multicolumn{6}{l}{HGA proposed by \cite{khameneh2023non}}                                                                                                                                                                                           &  & \multicolumn{6}{l}{IWO proposed by \cite{khameneh2023non}}                                                                                                                                                                                            \\ 
\cline{6-9}\cline{11-16}\cline{18-23}\cline{25-30}
                       &                        &                       &                        &  & Obj             & UB               & \begin{tabular}[c]{@{}l@{}}Time\\(sec)\end{tabular} & \begin{tabular}[c]{@{}l@{}}Gap\\(\%)\end{tabular} &  & Obj             & \begin{tabular}[c]{@{}l@{}}Time\\(sec)\end{tabular} & \begin{tabular}[c]{@{}l@{}}Imp\\(\%)\end{tabular} & CV              & SR            & \begin{tabular}[c]{@{}l@{}}DG\\(\%)\end{tabular} &  & Obj             & \begin{tabular}[c]{@{}l@{}}Time\\(sec)\end{tabular} & \begin{tabular}[c]{@{}l@{}}Imp\\(\%)\end{tabular} & CV              & SR            & \begin{tabular}[c]{@{}l@{}}DG\\(\%)\end{tabular} &  & Obj             & \begin{tabular}[c]{@{}l@{}}Time\\(sec)\end{tabular} & \begin{tabular}[c]{@{}l@{}}Imp\\(\%)\end{tabular} & CV              & SR            & \begin{tabular}[c]{@{}l@{}}DG\\(\%)\end{tabular}  \\ 
\hline
C1                     & 1                      & 3                     & 1                      &  & 28              & 28               & 13                                                  & 0.00                                              &  & 28.0            & 0.58                                                & 0.00                                              & 0.0000          & 1.0           & 0.00                                                &  & 28.0            & 4.99                                                & 0.00                                              & 0.0000          & 1.0           & 0.00                                                &  & 28.0            & 0.06                                                & 0.00                                              & 0.0000          & 1.0           & 0.00                                                 \\
                       & 2                      & 5                     & 1                      &  & 40              & 40               & 223                                                 & 0.00                                              &  & 40.0            & 1.60                                                & 0.00                                              & 0.0000          & 1.0           & 8.11                                                &  & 40.0            & 6.45                                                & 0.00                                              & 0.0000          & 1.0           & 8.11                                                &  & 38.0            & 1.89                                                & -5.00                                             & 0.0392          & 0.5           & 8.11                                                 \\
                       & 3                      & 7                     & 2                      &  & 45              & 45               & 6618                                                & 0.00                                              &  & 45.0            & 2.21                                                & 0.00                                              & 0.0000          & 1.0           & 4.65                                                &  & 45.0            & 14.50                                               & 0.00                                              & 0.0000          & 1.0           & 4.65                                                &  & 44.0            & 4.09                                                & -2.22                                             & 0.0339          & 0.7           & 4.65                                                 \\
                       & 4                      & 9                     & 2                      &  & 60              & 60               & 1757                                                & 0.00                                              &  & 60.0            & 3.70                                                & 0.00                                              & 0.0000          & 1.0           & 0.00                                                &  & 60.0            & 15.28                                               & 0.00                                              & 0.0000          & 1.0           & 0.00                                                &  & 59.8            & 7.28                                                & -0.33                                             & 0.0071          & 1.0           & 0.00                                                 \\
                       & 5                      & 11                    & 2                      &  & 130             & 198              & 43200                                               & 52.31                                             &  & 132.8           & 16.21                                               & 2.15                                              & 0.0032          & 1.0           & 6.40                                                &  & 132.1           & 21.93                                               & 1.62                                              & 0.0083          & 1.0           & 6.40                                                &  & 128.3           & 7.77                                                & -1.31                                             & 0.0152          & 0.4           & 4.80                                                 \\
                       & 6                      & 13                    & 3                      &  & 109             & 250              & 43200                                               & 129.36                                            &  & 112.5           & 7.92                                                & 3.21                                              & 0.0096          & 0.9           & 8.49                                                &  & 110.8           & 34.72                                               & 1.65                                              & 0.0038          & 0.0           & 4.72                                                &  & 107.5           & 5.98                                                & -1.38                                             & 0.0133          & 0.0           & 3.77                                                 \\
                       & 7                      & 15                    & 3                      &  & 173             & 366              & 43200                                               & 111.56                                            &  & 193.1           & 32.12                                               & 11.62                                             & 0.0079          & 0.9           & 6.52                                                &  & 184.0           & 43.59                                               & 6.36                                              & 0.0120          & 0.0           & 1.63                                                &  & 177.2           & 12.81                                               & 2.43                                              & 0.0118          & 0.0           & -2.17                                                \\
                       & 8                      & 17                    & 3                      &  & 168             & 359              & 43200                                               & 113.69                                            &  & 180.5           & 35.66                                               & 7.44                                              & 0.0129          & 0.8           & 8.93                                                &  & 173.3           & 48.14                                               & 3.15                                              & 0.0116          & 0.0           & 4.76                                                &  & 166.9           & 15.29                                               & -0.65                                             & 0.0365          & 0.0           & 4.76                                                 \\
                       & 9                      & 19                    & 4                      &  & 271             & 507              & 43200                                               & 87.08                                             &  & 288.5           & 22.46                                               & 6.46                                              & 0.0041          & 1.0           & 5.84                                                &  & 271.0           & 61.92                                               & 0.00                                              & 0.0154          & 0.0           & 1.09                                                &  & 258.4           & 15.66                                               & -4.65                                             & 0.0167          & 0.0           & -2.55                                                \\
                       & 10                     & 21                    & 4                      &  & 268             & 571              & 43200                                               & 113.06                                            &  & 334.5           & 30.35                                               & 24.81                                             & 0.0133          & 0.8           & 6.29                                                &  & 325.6           & 67.63                                               & 21.49                                             & 0.0109          & 0.0           & 3.77                                                &  & 301.6           & 25.50                                               & 12.54                                             & 0.0244          & 0.0           & -0.63                                                \\ 
\hline
\textbf{Average}       &                        &                       &                        &  & \textbf{129.20} & \textbf{242.40}  & \textbf{26781.10}                                   & \textbf{60.71}                                    &  & \textbf{141.49} & \textbf{15.28}                                      & \textbf{5.57}                                     & \textbf{0.0051} & \textbf{0.94} & \textbf{5.52}                                       &  & \textbf{136.98} & \textbf{31.92}                                      & \textbf{3.43}                                     & \textbf{0.0062} & \textbf{0.50} & \textbf{3.51}                                       &  & \textbf{130.97} & \textbf{9.63}                                       & \textbf{-0.06}                                    & \textbf{0.0198} & \textbf{0.36} & \textbf{2.07}                                        \\ 
\hline
C2                     & 1                      & 23                    & 5                      &  & 269             & 530              & 43200                                               & 97.03                                             &  & 326.1           & 35.55                                               & 21.23                                             & 0.0062          & 1.0           & 7.49                                                &  & 316.4           & 81.56                                               & 17.62                                             & 0.0212          & 0.2           & 6.84                                                &  & 297.4           & 23.41                                               & 10.56                                             & 0.0310          & 0.0           & 1.63                                                 \\
                       & 2                      & 25                    & 5                      &  & 313             & 649              & 43200                                               & 107.35                                            &  & 375.7           & 47.82                                               & 20.03                                             & 0.0097          & 0.8           & 5.26                                                &  & 348.7           & 89.92                                               & 11.41                                             & 0.0104          & 0.0           & -1.11                                               &  & 335.1           & 39.93                                               & 7.06                                              & 0.0171          & 0.0           & -4.43                                                \\
                       & 3                      & 27                    & 6                      &  & 308             & 740              & 43200                                               & 140.26                                            &  & 518.1           & 64.90                                               & 68.21                                             & 0.0072          & 1.0           & 6.07                                                &  & 490.2           & 124.53                                              & 59.16                                             & 0.0190          & 0.0           & 1.82                                                &  & 445.6           & 37.52                                               & 44.68                                             & 0.0369          & 0.0           & -2.23                                                \\
                       & 4                      & 29                    & 6                      &  & 299             & 795              & 43200                                               & 165.89                                            &  & 490.0           & 75.03                                               & 63.88                                             & 0.0077          & 0.8           & 10.18                                               &  & 458.7           & 140.95                                              & 53.41                                             & 0.0190          & 0.0           & 3.54                                                &  & 431.8           & 41.15                                               & 44.41                                             & 0.0279          & 0.0           & -0.66                                                \\
                       & 5                      & 31                    & 6                      &  & 315             & 744              & 43200                                               & 136.19                                            &  & 439.5           & 70.02                                               & 39.52                                             & 0.0127          & 0.6           & 8.47                                                &  & 416.0           & 146.16                                              & 32.06                                             & 0.0150          & 0.0           & 3.15                                                &  & 391.5           & 57.50                                               & 24.29                                             & 0.0220          & 0.0           & -2.66                                                \\
                       & 6                      & 33                    & 7                      &  & 303             & 743              & 43200                                               & 145.22                                            &  & 508.6           & 53.72                                               & 67.85                                             & 0.0092          & 0.9           & 6.63                                                &  & 481.5           & 157.87                                              & 58.91                                             & 0.0186          & 0.0           & 4.14                                                &  & 450.7           & 44.43                                               & 48.75                                             & 0.0281          & 0.0           & -3.31                                                \\
                       & 7                      & 35                    & 7                      &  & 436             & 1000             & 43200                                               & 129.36                                            &  & 566.6           & 111.40                                              & 29.95                                             & 0.0177          & 0.5           & 6.24                                                &  & 547.2           & 179.75                                              & 25.50                                             & 0.0202          & 0.0           & 3.67                                                &  & 482.2           & 61.72                                               & 10.60                                             & 0.0420          & 0.0           & -5.14                                                \\
                       & 8                      & 37                    & 7                      &  & 497             & 855              & 43200                                               & 72.03                                             &  & 626.9           & 285.32                                              & 26.14                                             & 0.0083          & 0.8           & 4.78                                                &  & 608.6           & 188.29                                              & 22.45                                             & 0.0103          & 0.0           & 2.14                                                &  & 581.9           & 65.94                                               & 17.08                                             & 0.0163          & 0.0           & -1.65                                                \\
                       & 9                      & 39                    & 8                      &  & 370             & 821              & 43200                                               & 121.89                                            &  & 542.7           & 234.30                                              & 46.68                                             & 0.0048          & 1.0           & 5.80                                                &  & 527.8           & 244.56                                              & 42.65                                             & 0.0172          & 0.0           & 3.48                                                &  & 487.9           & 66.40                                               & 31.86                                             & 0.0196          & 0.0           & -3.09                                                \\
                       & 10                     & 41                    & 8                      &  & 513             & 1062             & 43200                                               & 107.02                                            &  & 691.8           & 195.60                                              & 34.85                                             & 0.0093          & 0.9           & 4.00                                                &  & 643.2           & 206.64                                              & 25.38                                             & 0.0142          & 0.0           & -3.26                                               &  & 617.3           & 103.54                                              & 20.33                                             & 0.0125          & 0.0           & -6.96                                                \\ 
\hline
\textbf{Average}       &                        &                       &                        &  & \textbf{362.30} & \textbf{793.90}  & \textbf{43200.00}                                   & \textbf{122.22}                                   &  & \textbf{508.60} & \textbf{117.37}                                     & \textbf{41.83}                                    & \textbf{0.0093} & \textbf{0.83} & \textbf{6.49}                                       &  & \textbf{483.83} & \textbf{156.02}                                     & \textbf{34.86}                                    & \textbf{0.0165} & \textbf{0.02} & \textbf{2.44}                                       &  & \textbf{452.14} & \textbf{54.15}                                      & \textbf{25.96}                                    & \textbf{0.0253} & \textbf{0.00} & \textbf{-2.85}                                       \\ 
\hline
C3                     & 1                      & 43                    & 9                      &  & 310             & 1038             & 43200                                               & 234.84                                            &  & 610.4           & 123.41                                              & 96.90                                             & 0.0076          & 0.8           & 5.46                                                &  & 581.0           & 251.47                                              & 87.42                                             & 0.0101          & 0.0           & 0.85                                                &  & 530.6           & 58.09                                               & 71.16                                             & 0.0323          & 0.0           & -5.63                                                \\
                       & 2                      & 45                    & 9                      &  & 254             & 1042             & 43200                                               & 310.24                                            &  & 671.7           & 180.09                                              & 164.45                                            & 0.0115          & 0.6           & 6.38                                                &  & 644.0           & 287.41                                              & 153.54                                            & 0.0128          & 0.0           & 1.87                                                &  & 616.5           & 118.96                                              & 142.72                                            & 0.0133          & 0.0           & -1.40                                                \\
                       & 3                      & 47                    & 10                     &  & 287             & 980              & 43200                                               & 241.46                                            &  & 673.1           & 321.96                                              & 134.53                                            & 0.0086          & 0.9           & 5.76                                                &  & 655.5           & 343.03                                              & 128.40                                            & 0.0114          & 0.0           & 3.12                                                &  & 620.9           & 141.32                                              & 116.34                                            & 0.0130          & 0.0           & -1.25                                                \\
                       & 4                      & 49                    & 10                     &  & 271             & 1248             & 43200                                               & 360.52                                            &  & 906.6           & 820.88                                              & 234.54                                            & 0.0069          & 0.9           & 3.38                                                &  & 862.1           & 358.23                                              & 218.12                                            & 0.0104          & 0.0           & -1.13                                               &  & 822.1           & 108.91                                              & 203.36                                            & 0.0172          & 0.0           & -4.17                                                \\
                       & 5                      & 51                    & 10                     &  & 526             & 1314             & 43200                                               & 149.81                                            &  & 884.9           & 755.18                                              & 68.23                                             & 0.0058          & 1.0           & 4.43                                                &  & 849.1           & 370.12                                              & 61.43                                             & 0.0118          & 0.0           & 0.47                                                &  & 804.8           & 132.88                                              & 53.00                                             & 0.0188          & 0.0           & -3.03                                                \\
                       & 6                      & 53                    & 11                     &  & NS\textsuperscript{1}               & 1356             & 43200                                               & -                                                 &  & 843.5           & 1049.92                                             & -                                                 & 0.0105          & 0.7           & 6.46                                                &  & 784.7           & 427.21                                              & -                                                 & 0.0118          & 0.0           & -0.50                                               &  & 752.9           & 144.12                                              & -                                                 & 0.0139          & 0.0           & -3.48                                                \\
                       & 7                      & 55                    & 11                     &  & 324             & 1464             & 43200                                               & 351.85                                            &  & 963.1           & 823.62                                              & 197.25                                            & 0.0066          & 0.9           & 2.65                                                &  & 929.5           & 426.51                                              & 186.88                                            & 0.0098          & 0.0           & 0.32                                                &  & 870.2           & 114.62                                              & 168.58                                            & 0.0090          & 0.0           & -6.78                                                \\
                       & 8                      & 57                    & 11                     &  & NS              & 1255             & 43200                                               & -                                                 &  & 837.3           & 770.33                                              & -                                                 & 0.0063          & 1.0           & 4.44                                                &  & 803.2           & 441.74                                              & -                                                 & 0.0147          & 0.0           & 1.36                                                &  & 758.3           & 158.05                                              & -                                                 & 0.0163          & 0.0           & -4.32                                                \\
                       & 9                      & 59                    & 12                     &  & 410             & 1481             & 43200                                               & 261.22                                            &  & 965.9           & 1026.06                                             & 135.59                                            & 0.0051          & 1.0           & 5.64                                                &  & 925.7           & 484.51                                              & 125.78                                            & 0.0122          & 0.0           & 1.84                                                &  & 874.1           & 142.64                                              & 113.20                                            & 0.0114          & 0.0           & -3.90                                                \\
                       & 10                     & 61                    & 12                     &  & NS              & 1630             & 43200                                               & -                                                 &  & 1017.4          & 757.15                                              & -                                                 & 0.0088          & 0.7           & 3.50                                                &  & 961.7           & 499.38                                              & -                                                 & 0.0144          & 0.0           & -0.70                                               &  & 905.1           & 177.66                                              & -                                                 & 0.0206          & 0.0           & -5.70                                                \\ 
\hline
\textbf{Average}       &                        &                       &                        &  & \textbf{340.29} & \textbf{1223.86} & \textbf{43200.00}                                   & \textbf{272.85}                                   &  & \textbf{837.39} & \textbf{662.86}                                     & \textbf{147.36}                                   & \textbf{0.0078} & \textbf{0.85} & \textbf{4.81}                                       &  & \textbf{799.65} & \textbf{388.96}                                     & \textbf{137.37}                                   & \textbf{0.0119} & \textbf{0.00} & \textbf{0.75}                                       &  & \textbf{755.55} & \textbf{129.73}                                     & \textbf{124.05}                                   & \textbf{0.0166} & \textbf{0.00} & \textbf{-3.97}                                       \\ 
\hline
\multicolumn{30}{l}{\textsuperscript{1}NS: No feasible solution found in 43200 seconds}                                                                                                                                                                                                                                                                                                                                                                                                                                                                                                                                                                                                                                                                                                                                                                                                                                                                                 
\end{tabular}}

\end{minipage}
}

\end{table}


%% file: Tables/Table_SensAnalLevels.tex
\begin{table}[h]
\centering
\caption{Parameter levels used in the sensitivity analysis.}
\label{Table_SensitivityLevels}
\def\arraystretch{1.3}
\scalebox{0.8}{
\begin{tblr}{
  colspec = {l p{2cm} p{2cm} p{2cm}},
  row{2} = {c},
  cell{1}{1} = {r=2}{},
  cell{1}{2} = {c=3}{c},
  cell{3}{2} = {c},
  cell{3}{3} = {c},
  cell{3}{4} = {c},
  cell{4}{2} = {c},
  cell{4}{3} = {c},
  cell{4}{4} = {c},
  cell{5}{2} = {c},
  cell{5}{3} = {c},
  cell{5}{4} = {c},
  cell{6}{1} = {c=4}{},
  hline{1,3,6} = {-}{},
  hline{2} = {2-4}{},
}
Parameters                                                & Levels &        &      \\
                                                          & Low    & Medium & High \\
Payload capacity (number of blood bags)                   & 5      & 10     & 15   \\
Battery endurance (minutes)                               & 15     & 30     & 45   \\
Speed\textsuperscript{1}                                  & 1      & 2      & 3    \\

{\textsuperscript{1}speed is given as a multiplier relative to truck speed.}    
\end{tblr}
}
\end{table}